\documentclass[11pt,onecolumn]{article}
\usepackage{amscd}
\usepackage{times}
\usepackage{amsmath}
\usepackage{amssymb}
\usepackage{xspace}
\usepackage{theorem}
\usepackage{graphicx}
\usepackage{ifpdf}
\usepackage{url,hyperref}
\usepackage{latexsym}
\usepackage{euscript}
\usepackage{xspace}
\usepackage[all]{xy}
\usepackage{color}
\usepackage{makeidx}
\usepackage{tikz}

\long\def\remove#1{}

\newtheorem{theorem}{Theorem}[section] 
\newtheorem{lemma}[theorem]{Lemma}

\newtheorem{obs}[theorem]{Observation}

\newtheorem{definition}[theorem]{Definition}
\newtheorem{proposition}[theorem]{Proposition}
\newenvironment{proof}{{\em Proof:}}{\hfill{\hfill\rule{2mm}{2mm}}}

\newcommand {\mm}[1] {\ifmmode{#1}\else{\mbox{\(#1\)}}\fi}

\newcommand{\B}                        {\mathrm {\mathbb{B}}}

\newcommand{\A}                        {\mathrm {\mathbb{A}}}

\newcommand{\img}{\mathrm {img}}
\newcommand{\dom}{\mathrm {dom}}
\newcommand{\reg}{\mathrm{reg}}

\newcommand{\mul}{\mathrm{mul}}

\newcommand{\cancel}[1]

\begin{document}

\title{Linear relations, monodromy 
and Jordan  cells of a circle valued map.}



\author{
Dan Burghelea  \thanks{
Department of Mathematics,
The Ohio State University, Columbus, OH 43210,USA.
Email: {\tt burghele@math.ohio-state.edu}}
}
\date{}
\maketitle

\begin{abstract}
\vskip .2in

In this paper we consider the definition of {\it monodromy of an angle valued map}
based on linear relations as  proposed in \cite{BH}. This definition  
provides an alternative  treatment of monodromy and computationally  an alternative calculation of the {\it Jordan cells}, 
topological persistence invariants of a circle valued maps introduced in  \cite {BD11}.

We give a new geometric proof  that the monodromy is actually a homotopy invariant of a pair  $(X,\xi)$
consisting of a compact ANR  $X$ and  an integral cohomology class $\xi\in H^1(X; \mathbb Z),$ without any reference to the infinite cyclic cover associated to $\xi,$ as in \cite {BH}, or to the graph representation associated an angle valued map defining  $\xi$ as in \cite{BD11}. 

Most important, we describe an algorithm to calculate the monodromy for a simplicial angle valued map defined on a finite simplicial 
complex, providing a new algorithm for the calculation of the Jordan cells of the map 
shorter than the one proposed in \cite {BD11}.

We  indicate the computational usefulness of {\it Jordan cells}, and in particular  of the proposed algorithm, for the calculation of other basic topological invariants of $(X,\xi).$ 
\end{abstract}

\thispagestyle{empty}
\setcounter{page}{1}

\tableofcontents



\section{Introduction}

Absolute neighborhood retracts (ANRs) are topological spaces $X$  which whenever  $ i:X\to Y$ is an embedding into a normal topological space $Y$ there exists a neighborhood $U$ of $i(X)$ in $Y$ and a retraction of $U$ onto $i(X)$ cf \cite {Hu}.
\vskip .1in

Let $X$ be a compact ANR 
\footnote {The reader unfamiliar with the notion of ANR should always think to the  main examples, spaces homeomorphic to simplicial complexes or to a CW complexes},  
$\xi\in H^1(X; \mathbb Z)$  and $\kappa$ a field with algebraic closure $\overline \kappa.$ 
The $r-$monodromy, $r\in \mathbb Z_{\geq 0},$ is a 
similarity (= conjugacy)  class of linear isomorphism  $T^{(X,\xi)}(r): V_r(X,\xi) \to V_r(X;\xi),$ with $V_r(X,\xi)$ finite dimensional 
$\kappa-$vector spaces, cf definition \ref {SS21} below. The Jordan decomposition of a square matrix permits to assign to the linear isomorphism  $T^{(X,\xi)}(r)$ the collection  $\mathcal J_r(X;\xi)$ of pairs $(\lambda,k),$  $\lambda \in \overline \kappa\setminus 0, k\in \mathbb Z_{\geq 1}$, referred to as {\it Jordan cells} in dimension $r.$ They provide a complete sets of invariants of the similarity class of $T^{(X,\xi)}(r).$ 
If $f:X\to \mathbb S^1$ is a tame map as in \cite {BD11} and $\xi_f$ the cohomology class defined by $f,$ then the set $\mathcal  J_r(X;\xi_f)$ coincides with  the set of Jordan cells $\mathcal J_r(f)$ considered in \cite{BD11} in relation with 
the  topological persistence of  the  circle valued map $f,$ cf  \cite{BH}.
 \vskip .1in   
Recall that topological persistence for a real or circle valued map $f:X\to \mathbb R$ or $f:X \to \mathbb S^1$  ($\mathbb S^1= \mathbb R/2\pi \mathbb Z$)   analyses the changes  in the homology of the levels $f^{-1}(\theta), \theta \in \mathbb R \ \mbox {or}\  \mathbb S^1.$  It records the {\it detectability} and the {\it death}  of  homology of the levels  in terms of {\it bar codes} cf  \cite {BD11} or \cite {CSD09}. In case of a circle valued map, in addition to {\it death} and {\it detectability},  there is an additional feature of interest to be recorded, the {\it return of some homology classes}  of $f^{-1}(\theta)$ when the angle $\theta $ increases or decreases with  $2\pi.$ 
This feature is recorded as {\it Jordan cells} which were introduced in \cite {BD11} and describe what the topologists refer to  as the {\it homological monodromy} or simply the {\it monodromy}. 
In \cite{BH} we have proposed a definition for  {\it monodromy} and implicitly for {\it Jordan cells}  based on {\it linear relations}. For the purpose of this paper  the needed backgrund on {\it Linear Relations} is presented in section \ref{S3}; for more background  the reader can consult \cite {SSW05} or \cite{BH} section 8.
\vskip .1in

 In this paper we review the definition of monodromy based on linear relations, provide a new geometric proof of its homotopy invariance (without any reference to Novikov homology used in \cite{BH}) and, more important, propose a new algorithm  for the calculation of $\mathcal J_r(f)=\mathcal J_r(X;\xi_f),$ for $X$ a finite  simplicial complex and $f$ a simplicial map.  The notation $\mathcal J_r(X;\xi_f)$ indicates that the collection of pairs $(\lambda, k)$ depend on the pair $(X,\xi\in H^1(X;\mathbb Z))$ rather than $f:X\to \mathbb S^1$ representing $\xi.$ 
 
In the present  approach the monodromy is first defined for a  continuous map $f:X\to \mathbb S^1$ and a {\it weakly regular angle} $\theta\in \mathbb S^1$ (see the definitions in section 3). Note that not all compact ANR's  admit  angle valued maps with weakly regular angles, cf \cite{DW}.
Note also that for a simplicial map all angles are weakly regular. We reduce the general case of an arbitrary  compact ANR and a continuous map,  
which might have no weakly regular angles, to the case of simplicial complexes and simplicial maps based on results on the topology of compact Hilbert cube manifolds.

Proposition \ref{T31} shows that the monodromy proposed is independent of the weakly regular angle, remains the same for 
maps  which have weakly regular angles and are homotopic and does not change when one replaces the map  by its composition with the projection $X\times K\to X,$  $K$ an acyclic compact ANR. These facts ultimately show that the monodromy can be associated to a pair $(X,\xi \in H^1(X;\mathbb Z)),$  $X$ any compact ANR, and the assignment is a homotopy invariant of the par $(X,\xi),$  cf Theorem \ref{T36}.  All these facts are established in section 3, based on elementary linear algebra of linear relations summarized in section 2. They also follow from the definition of monodromy based on the homology of the infinite cyclic cover associated to $\xi$ described in \cite {BH} but this is exactly what the present treatment wants to avoid.

The algorithm for calculating $\mathcal J_r(f)$ for $f$ a simplicial angle valued map is discussed in section \ref{S4}. 

In section \ref{S5} we  indicate applications of the calculation of Jordan cells to the calculation topological invariants whose standard definitions involve infinite cyclic covers, {\it computer unfriendly} objects. 
In particular one provides new ways to calculate  Novikov Betti numbers, the Alexander polynomial of a knot and a few other  invariants, see section \ref{S5}. 

In section 3 we notice that  a generalization of the homological  monodromy discussed in this paper can be obtained when the  singular homology $H_r$ is replaced by a vector space valued homotopy functor $F$ which is half exact in the sense of A. Dold cf \cite {D}. This $F-$ monodromy is not investigated in this paper  but it might deserve attention \footnote{  
A slightly more general situation, when $X$ is equipped with a flat
bundle of $\mathcal A-$modules $W$,  $\mathcal A$ a finite type von Neunmann algebra, and $W$ a finite type Hilbert module will be considered in a sequel of this paper in the special case $\mathcal A$ is the von Neumann algebra $\mathcal N(\pi_1(M))$ and $W= \l^2(\pi_1(M))\otimes \mathbb C^n).$  Such monodromy will be used to the description of $L_2-$torsion in geometrically interesting situations.}

{\it Acknowledgements:} 
The idea of describing the Jordan cells considered in \cite {BD11} using linear relations  belongs to  Stefan Haller and was pursued in \cite{BH} not yet in print.

It is a pleasure to thank  S.Ferry for help in relation with the Appendix 2. and for bringing to our attention the reference \cite {DW}.

\section {Linear relations}\label {S1}

Fix a field $\kappa$ and let  $\tilde \kappa$ be its algebraic closure. 

\subsection{\bf Generalities on linear relations} \label {SS21}

Recall from \cite {SSW05} and \cite{BH}:

--  A linear relation $R:V_1\rightsquigarrow V_2$ is a linear subspace $R\subseteq V_1\times V_2.$ One writes $v_1 R   v_2$  iff $(v_1, v_2)\in R,$ $v_i\in V_i.$

{\it Examples: }
\begin{enumerate}
\item  Two  linear maps $\xymatrix{V_1\ar[r]^\alpha &W &V_2\ar[l]_\beta}$ provide the relation 

\hskip .5in  $R(\alpha,\beta) \subset V_1\times V_2 := \{(v_1, v_2) \mid \alpha(v_1)= \beta(v_2) \}.$
\item Two linear maps $\xymatrix {V_1 & U\ar[l]_a \ar[r]^b&V_2}$ provide the relation 

\hskip .5in $R<a,b> \subset V_1\times V_2 := \{(v_1, v_2) \mid  \exists u, a(u)= v_1, b(u)= v_2\}.$
\end{enumerate} 

--Two liner relations $R_1\colon V_1\rightsquigarrow V_2$  and $R_2\colon V_2\rightsquigarrow V_2$ can be composed in an obvious way,  $(v_1 (R_2\cdot R_1) v_3$ iff $\exists v_2$ such that $v_1 R_1 v_2$ and $v_2 R_2 v_3.$ The diagonal $\Delta \subset V\times V$ is playing the role of the identity.

--Given a linear relation $R: V_1\rightsquigarrow V_2$ denote by $R^\dag :V_2\rightsquigarrow V_1$ the relation defined by the property 
$v_2R^\dag v_1$ iff $v_1Rv_2.$ 
Clearly $(R_1\cdot R_2)^\dag= R^\dag_2\cdot R^\dag_1$ and $R^{\dag \dag} = R.$ 

The familiar category of finite dimensional vector spaces and linear maps can be extended to incorporate all linear relations as morphisms. The linear map $f:V_1\to V_2$ can be interpreted as the relation {\it $ \it{graph}\  f  \subset V_1\times V_2$} denoted by $\boxed{R(f)}= R(f,id_{V_2}),$  providing the embedding of the category of vector spaces and linear maps in the category of vector spaces and linear relations.  This extended category remains abelian.

--The direct sums $R' \oplus R'':V'_1\oplus V''_1\rightsquigarrow V_2'\oplus V''_2$ of two relations $R' :V'_1\rightsquigarrow V_2'$ and $R'':V_1'' \rightsquigarrow V_2''$  is defined in the obvious way, $(v'_1,v''_1) (R'\oplus R'') (v'_2, v''_2)$ iff $(v'_1 R' v'_2) $ and $(v''_1 R'' v''_2).$ 

--The relation with the same source and target $R' \colon V' \rightsquigarrow V'$ and $R'' \colon V'' \rightsquigarrow V''$ are similar  and one writes 
$\boxed {R'\sim R''}$ if  there exists the linear isomorphisms $\alpha: V'\to V''$  s.t. $R''\cdot R(\alpha)= R(\alpha)\cdot R'.$

Recall that two linear endomorphisms  $T:V\to V$ and $T':V'\to V'$ are called similar if there exists a linear isomorphism $C:V\to V'$ s.t. $C^{-1}\cdot T'\cdot C= T$.  One writes $T\sim T'$  if $T$ and $T'$ are similar and one denotes the similarity class 
of $T:V\to V$ by $[T];$  so $T\sim T'$ and $[T]= [T']$ mean the same thing.

As in the case of linear maps one denotes the similarity class
of the relation $R\colon V\rightsquigarrow V$ by $[R].$  Clearly   
when $T:V\to V$ is a linear map both notations $[T]$ and $[R(T)]$ mean  the same thing.
\vskip .1in

\noindent A linear relation $R\colon V\rightsquigarrow W$ gives rise to the following subspaces:
\begin{align*}
\dom(R)&:=\{v\in V\mid\exists w\in W:vRw\} = pr_V(R)
\\
\img(R)&:=\{w\in W\mid\exists v\in V:vRw\} = pr_W(R)  
\\
\ker(R)&:=\{v\in V\mid vR0\} {\cong} V\times 0 \cap R
\\
\mul(R)&:=\{w\in W\mid 0Rw\} {\cong} 0\times W \cap R
\end{align*}
Here   $pr_V$ and $pr_W$ denote the projections of $V\times W$ on $V$  and $W.$ 
We have 
\begin{obs}\
\begin{enumerate}
\item 
$\ker(R)\subseteq \dom(R)\subseteq V $ 
and
$W\supseteq \img(R)\supseteq \mul(R),$
\item $\ker(R^\dag)= \img (R)$ 
and
$\dom(R^\dag)= \img(R),$
\item 
$\dim\dom(R)+\dim\ker(R^\dag)=\dim(R)=\dim(R^\dag)= \dim\dom(R^\dag)+\dim\ker(R).$
\end{enumerate}
\end{obs}

It is  immediate, in view of the above definitions and above observation that :

\begin{lemma}\label {L22}\ 

1. A linear relation $R\colon V\rightsquigarrow W$ is of the form $R(f)$ for  $f:V\to W$ linear map iff 
$\dom R= V$ and $\mul R=0.$

2. A  linear relation $R\colon V \rightsquigarrow V$  is of the form $R(T)$ for $T:V\to V$ a linear isomorphism iff 
$\dom R= V$ and $\ker R=0.$ 
\end{lemma}

Let $R: V\rightsquigarrow V$ be a linear relation. 
Define
\begin{enumerate}
\item $D:\{v \in V  \mid  \exists v_i\in V,  i\in \mathbb Z, v_iR v_{i+1}, v_0= v\}.$ The relation  $R$ restricts to a relation $R_D :D\rightsquigarrow D.$
\item $K_+:=\{ v\in V \mid \exists v_i, i\in \mathbb Z_{\geq 0} , v_iR v_{i+1}, v_0= v\}.$
\item $K_+:=\{ v\in V \mid \exists v_i, i\in \mathbb Z_{\geq 0} , v_iR v_{i+1}, v_0= v\}.$
\item $V_{reg}:= \frac{D} {D \cap (K_+ + K_-)},\ $  $\pi: D\to \frac{D} {D \cap (K_+ + K_-)}$ the quotient map and $\iota: D\to V$ the inclusion.
\end{enumerate} 

Consider the composition of relations $$R_D= R(\iota)^\dag \cdot R\cdot R(\iota)$$ 
and define $$R_{reg}:= R(\pi) \cdot R_D \cdot R(\pi)^\dag: V_{reg}\rightsquigarrow V_{reg}.$$ 

\begin{proposition}\label{P23} (cf \cite{BH}) \ 
\begin{enumerate}
\item  There exists a linear isomorphism  $T^R: V_{reg}\to V_{reg}$ such that $R_{reg}= R(T^R).$

\item If $R\colon V\rightsquigarrow V$ and $R'\colon V'\rightsquigarrow V'$ are similar relations, i.e. there exists an isomorphism of vector spaces $\omega:V\to V'$ such that $R'= R(\omega)\cdot R \cdot R(\omega^{-1}),$  then $T^R$ and $T^{R'}$  are similar linear isomorphisms (i.e. 
$T^{R'}=\underline \omega \cdot T^R\cdot {\underline \omega}^{-1}$ for some isomorphism $\underline \omega$ ).
\item $R_{reg}^{-1}= (R^\dag)_{reg}.$

\item  $(R'\oplus R'')_{reg}= R'_{reg} \oplus R''_{reg}.$
\item Suppose $R_i: V_i\rightsquigarrow V_{i+1}, i=1,2, \cdots k$ with $V_1= V_{k+1}$ then 
$(R_i\cdots R_{i-1}\cdots R_1\cdot R_k\cdot R_{k-1}\cdots R_{i+1})_{reg}\sim (R_{k} \cdot R_{k-1}\cdots R_2\cdot R_1)_{reg}.$
\end{enumerate}
\end{proposition}

In view of the definition of $R_\reg$ it is immediate that :
\begin{obs}\
\begin{enumerate}
\item If $\alpha, \beta:V\to W$ are two isomorphisms then $T^{R(\alpha, \beta)}= \beta^{-1}\cdot \alpha.$
\item If $f:V\to V$ is a linear map and  $V_0$ is the generalized eigen-space of the eigenvalue $0$ then:

 $f(V_0)\subset V_0$,

$f$ induces $\hat f:V/V_0\to V/V_0$ and 

$T^{R(f)}\sim \hat f:V/V_0\to V/V_0.$ 
 \end{enumerate}
 \end{obs}

The following  technical Proposition will be used in  section \ref{SS42},  where an algorithm for  the calculation of $R(a,b)_{reg},$   
part of an algorithm for the calculation of the $r-$ monodromy will be presented. 

\begin{proposition}\label {P24}\
\begin{enumerate}
\item 
Consider the diagram 
\begin{equation}
\xymatrix { 
V\ar[r]^\alpha& W & V\ar[l]_\beta\\
V'\ar[r]^{\alpha'}\ar[u]^{\subseteq}& W' \ar[u]^{\subseteq}& V'\ar[l]_{\beta'}\ar[u]^{\subseteq}}
\end{equation}
and suppose that: 
 
$W'\supseteq \img{\alpha}\cap  \img{\beta}$

$V' = \alpha^{-1}(W')\cap \beta^{-1}(W')$ and 

$\alpha'$ and $\beta'$ the restriction of $\alpha$ and $\beta.$

Then $R(\alpha, \beta)_{reg}= R(\alpha', \beta')_{reg}.$

\item Consider the diagram 
\begin{equation}
\xymatrix { 
V\ar[r]^\alpha \ar[d]^{p'}& W \ar[d]^p& V\ar[l]_\beta \ar[d]^{p'}\\
V' \ar[r]^{\alpha'} & W'& V'\ar[l]_{\beta'} }
\end{equation}

with both $\alpha$  and $\beta$  surjective. Define :

$V'= V/\ker \alpha$,  $W'= W/ \beta(\ker\alpha),$

$p: W\to W' ,$ $p': V\to V'$ the canonical quotient maps,

$\overline \alpha : V' \to W$ induced from $\alpha,$ $\alpha'= p \cdot \overline \alpha ,$ 

$\beta'$ induced by passing to quotient from $\beta.$ 

Then $R(\alpha, \beta)_{reg}= R(\alpha', \beta')_{reg}.$
\end{enumerate}
\end{proposition}
For the reader's convenience the proofs of Propositions \ref{P23} and \ref{P24} 
are included in Appendix 1. 

\subsection {\bf Jordan cells,  characteristic  polynomial and the characteristic divisors}\label {sec:2.2}

Recall that a Jordan matrix $T(\lambda,k)$ is determined by a pair $(\lambda, k),$ $\lambda\in \overline \kappa$ and $k$ a positive integer.
When  $\lambda\ne 0$  the pair $(\lambda, k)$ is  called  in \cite {BD11}  {\it Jordan cell}. 
\begin{equation*}
T(\lambda; k)=
\begin{pmatrix}
\lambda & 1       & 0      & \cdots  & 0      \\
0       & \lambda & 1      & \ddots  & \vdots \\
0       & 0       & \ddots & \ddots  & 0      \\
\vdots  & \ddots  & \ddots & \lambda & 1      \\
0       & \cdots  & 0      & 0       & \lambda 
\end{pmatrix} _.
\end{equation*}
Any invertible square $n\times n-$ matrix is conjugated with a direct sum of Jordan cells (by Jordan decomposition theorem , cf \cite {G}) 
with $\lambda$  eigenvalue of the matrix. In different words  any conjugacy class of linear isomorphisms $T:V\to V,$ denoted by $[T],$ is determined by a unique collection of pairs, the Jordan cells $\mathcal J([T])$ or $\mathcal J(T).$  
Note that Any such collection determines and is determined by the collection of monic polynomials 

$$P^T(z) | P^T_1(z) |P^T_2(z)|\cdots P^T_{n-1}(z)$$
 where $P^T(z)= \det (zI- T)$ and $P^T_i(z)$ is the greatest common divisor of all $(n-i)\times (n-i)-$ minors of $zI- T,$ cf \cite{G}. 
The polynomials $P^T(z) | P^T_1(z) |P^T_2(z)|\cdots P^T_{n-1}(z)$ \footnote {$P(z) | Q(z)$ means that there exists a polynomial $R(z)$ such that $P(z)= Q(z) R(z).$} do not involve the algebraic closure $\overline \kappa$ .
The precise relation between them and the elements of $\mathcal J([T])$ is given  in \cite {G}.
\begin {definition} The Jordan cells of the linear relation $R\colon V\rightsquigarrow V$ is the collection $\mathcal J([T^{R_\reg}]).$
\end{definition}  
\vskip .2in 

\section {Monodromy} \label{S3}

In this section the homology of a space $X$ is the singular homology with coefficients in a field $\kappa$ fixed once for all and is denoted by $H_r(X),$  $r=0,1,2,\cdots.$ 

 An {\it angle} is a complex number $\theta= e^{it} \in \mathbb C, t\in \mathbb R$ and the set of all angles is denoted  by $\mathbb S^1=\{\theta= e^{it}\mid t\in \mathbb R\}.$  The space of angles identified to  $\mathbb S^1,$  is equipped with the distance $$d(\theta_2, \theta_2)=\inf \{ |t_2- t_1|\mid
  e^{it_1}=\theta_1, e^{it_2}= \theta_2 \}.$$ 

 In this paper {\bf all real valued or angle valued maps  $f:X\to \mathbb R$ or $f:X\to \mathbb S^1$ are proper continuous  maps with $X$ an ANR.}  The properness of $f$  forces the space $X$ to be locally compact in the first case and compact in the second.    

-- A value $t\in \mathbb R$ or $\theta\in \mathbb S^1$  is {\it weakly regular}  if $f^{-1}(t)$ or  $f^{-1}(\theta)$ is an ANR, hence a compact ANR\footnote {A compact ANR has the homotopy type of finite simplicial complex.}.

-- A  map $f$ whose set of weakly regular values is not empty is called {\it good} and  a map with all values weakly regular is called {\it weakly tame}. For $X$ a (locally finite) simplicial complex any $\mathbb R$ or $\mathbb S^1-$ valued simplicial map $f$ is weakly tame. 

--  An  ANR \  $X$ whose set of  weakly tame  maps is dense in the space of all  maps with the $C^0-$fine 
topology    
\footnote{For this paper the concepts of {\it  good map, tame map and good ANR} will be considered  under the hypothesis that the space is  compact, in which case $C^0-$fine topology is the same as the familiar compact open topology} is called a {\it good ANR}. 
There exist  compact ANR's (actually compact homological n-manifolds, cf \cite {DW}) with no co-dimension one subsets which are ANR's, hence compact ANR's which are not {\it good }ANR's. 
The spaces  homeomorphic to  simplicial complexes, finite  dimensional topological manifolds, or Hilbert cube manifolds  (see Appendix 2 for definitions)  are all {\it good}  ANR's. The first because any continuous map can be approximated by simplicial 
 maps w.r. to a convenient subdivision, the last by more subtle reasons explained in Appendix 2.  
\vskip .1in 
 
 As pointed out in introduction,  the $r-${\it monodromy}, cf Definition \ref {D2} below,  will be first defined for good maps and will involve an angle $\theta,$ which is a weakly regular value. It will be shown that different choices  of such angles lead to the same $r-$monodromy, and that the $r-$monodromy  depends  only on the cohomology class $\xi_f$ associated with the map $f.$  
  
Once some elementary properties will be established 
for good ANRs and maps\footnote {actually it suffices to established them for simp;laical complexes ad simplicial maps}, using results on Hilbert cube manifolds, it will be shown that the $r-$monodromy can be associated to any angle valued map and is a homotopy invariants for any pair $(X,\xi\in H^1(X;\mathbb Z)),$ $X$ any  compact ANR. 

The following observations will be useful.

 \begin {proposition} \label{O52}\
 \begin{enumerate}
 \item  Two  maps $f,g:X\to \mathbb S^1$ with $D(f,g)= \sup_{x\in X} d(f(x), g(x))<\pi$ are  homotopic by a canonical homotopy, the  ''convex combination`` homotopy.
 \item  Suppose $X$ is a good ANR,   $f,g : X\to \mathbb S^1$ two homotopic angle valued maps and $\epsilon >0.$ Then there exists a finite collection of maps $f_0, f_1, \cdots f_k, f_{k+1},$ such that:
 
  a) $f_0= f, f_{k+1}=g,$
  
  b) $f_i$ are weakly tame maps  for $i=1,2,\cdots k,$
  
  c) $D(f_i, f_{i+1}) <\epsilon.$
  \end{enumerate}
\end{proposition}
 
Indeed if $f$ and $g$ are viewed as maps with values in $\mathbb C$ then the map $h_t(x)= \frac {t g(x)+(1-t) f(x)}{ |t g(x)+(1-t) f(x)|}$, $0\leq t\leq 1,$ provides the desired homotopy stated in item 1.  The condition $ D(f(x), g(x))<\pi$  insures that $|t g(x)+(1-t) f(x)|\ne 0.$  

Item 2. follows from the local contractibility of the space of maps when equipped with the distance $D.$

\subsection {\bf Real valued maps}\label {SS41}

For $f: X\to \mathbb R$ a  real valued map and 
$a\in \mathbb R$ denote by:

$X^f_a,$ the sub-level $X^f_a:= f^{-1}((-\infty, a])$; if $a$ is weakly regular value  then $X^f_a:= f^{-1}((-\infty, a])$is an ANR,

$X^a_f,$  the super-level $X_f^a:= f^{-1}([a,\infty);$ if $a$ is weakly regular  value then $X^f_a:= f^{-1}([a,\infty))$ is an ANR.

\noindent For $f: X\to \mathbb R$ and $g: X\to \mathbb R$   maps  $a<b$ s.t.  $f^{-1}(a)\subset g^{-1}(-\infty,b)$
denote by 
$$X^{f,g}_{a,b}:= X^g_b\cap X_f^a.$$  
If $b$ is a weakly regular value for $g$ and $a$ is  weakly regular value for $f$ then $X^{f,g}_{a,b}$ is a compact ANR. 
This insures that $H_r(g^{-1}(a)), H_r(f^{-1}(b))$ and  $H_r (X^{f,g}_{a,b})$ have finite dimension. 

Denote by  $R^{f,g}_{a,b}(r)$  the linear relation defined by the  linear maps  ${i_1 (r)}$ and ${i_2 (r)}$ induced by the inclusions 
$f^{-1}(a)\subset X^{f,g}_{a,b}$ and $g^{-1}(b)\subset  X^{f,g}_{a,b}.$
$$\xymatrix {H_r(f^{-1}(a)) \ar[r]^{i_1 (r)} & H_r(X^{f,g}_{a,b}) & H_r(g^{-1}(b))\ar[l]_{i_2 (r)}}.$$ 

\begin{proposition} \label {T30} Let $t_1< t_2 < t_3.$  Suppose that $t_1$ is weakly regular for $f,$  $t_2$ is  weakly regular for $g$ and $g^{-1}(t_2)\subset f^{-1}((t_1, t_3)).$ 
Then one has 
$$R^{g,f}_{t_2, t_3}(r) \cdot R^{f,g}_{t_1, t_2}= R^{f,f}_{t_1, t_3}(r).$$
\end{proposition}
\vskip .1in
\begin{proof} 
The verification is a consequence of the exactness of the following piece of 
 Meyer--Vietoris sequence 

\begin{equation}\label{E7}
\xymatrixcolsep{5pc}\xymatrix 
{H_r(g^{-1}(t_2))\ar[r]^{i'_1\oplus i'_2}&H_r(X^{f,g}_{t_1,t_2}) \oplus H_r(X^{g,f}_{t_2,t_3}) \ \ar[r]^{i_1 - i_2}& H_r(X^{f,f}_{t_1,t_3}}
\end{equation}
whose linear maps involved in the sequence (\ref {E7}) and in the commutative diagram  below are induced by obvious inclusions
\footnote {In order to  lighten the writing, in both (\ref{E7}) and (\ref{E8}), "$r$'' was dropped off  $i$'s and $i'$ 's the notations for  the inclusion induced linear maps in $r-$homology}.

\scriptsize{
\begin{equation}\label{E8}
\xymatrix{           &                                                  & H_r(X^{f,f}_{t_1, t_3}) &  &                                 \\
H_r(f^{-1}(t_1))\ar[r]^{j_1}\ar[rru]^{I_1} 
& H_r(X^{f,g}_{t_1, t_2})\ar[ru]_{i_1} & H_r(g^{-1}(t_2))\ar[l]^{i'_1}\ar[r]_{i'_2} & H_r(X^{g,f}_{t_2,t_3}) \ar[lu]^{i_2}&  H_r(f^{-1}(t_3))\ar[l]_{j_2}\ar[llu]_{I_2}} 
\end{equation}
} 
\normalsize
 \vskip .1in

Indeed 
the commutativity of the diagram (\ref{E8}) implies that 
$ x R^{f,f}_{t_1, t_2}y,$  for  $x\in H_r (f^{-1}(t_1))$ and $y\in H_r (f^{-1}(t_3))$ \  iff \  $i_1( j_1(x))- i_2 ( j_2(y))=0.$  

By the exactness of the sequence (\ref{E7})  one has \ \  
$i_1( j_1(x))- i_2 ( j_2(y))=0$  iff there exists $u\in H_r(g^{-1}(t_2))$ such that \ $(i'_1 \oplus i'_2)(u) = (j_1(x), j_2(y)).$  
This happens iff $x R^{f,g}_{t_1, t_2} u $ and $u R^{g,f}_{t_2, t_3} y.$ which means $x R^{f,f}_{t_1, t_2} y.$\ \
\end{proof}

\subsection {\bf Angle valued maps}\label {SS3.2}

Let $f: X\to \mathbb S^1$ be an angle valued map.  Let  
$u\in H^1(S^1;\mathbb Z)\equiv \mathbb Z$ be the generator defining the orientation  of $\mathbb S^1.$ Here $\mathbb S^1$ is regarded as an oriented one dimensional manifold.  Let 
$f^\ast: H^1(\mathbb S^1;\mathbb Z)\to H^1(X;\mathbb Z)$ be the  homomorphism induced by $f$ in integral cohomology  and $\xi_f= f^\ast(u) \in H^1(X;\mathbb Z).$  
  It is a well known fact in homotopy theory that the assignment $f\rightsquigarrow \xi_f$ establishes a bijective correspondence between  the set of homotopy classes of continuous maps from $X$ to $\mathbb S^1$ and $H^1(X;\mathbb Z).$
\vskip .2in

{\bf The cut at $\theta$ (with respect to  the map $f: X\to \mathbb S^1$) } 

\noindent For  $\theta \in \mathbb S^1,$  a weakly regular value for $f,$ 
 define 
 {\bf the cut at $\theta = e^{it}$}  to be the space $\overline X^f_\theta,$  the two sided compactification of $X\setminus f^{-1}(\theta)$ with  sides  $f ^{-1}(\theta).$ 
Precisely as a set $\overline X^f_\theta$  is a disjoint union of three parts,  
$\overline X^f_\theta = f^{-1}(\theta)(1) \sqcup f^{-1}(\mathbb S^1\setminus \theta) \sqcup f^{-1}(\theta)(2),$ with $f^{-1}(\theta)(1)$ and $f^{-1}(\theta)(2)$ two copies of $f^{-1}(\theta).$ 

The topology on $\overline X^f_\theta$ is the only topology  which makes $\overline X^f_\theta$ compact and the map from $\overline X^f_\theta$ to  $X$  defined by identity on each part  continuous and a homeomorphism on the image when restricted to each part. 
The compact space $\overline X^f_\theta$ is a compact ANR.  

The obvious inclusions $i_1, i_2,$ $\xymatrix{  f^{-1}(\theta )\ar[r]^{i_1} &\overline X_\theta  & f^{-1}(\theta) \ar[l]_{i_2}}$  induce in homology in dimension $r$ the  linear maps  (between finite dimensional vector spaces) $i_1(r)$ and $i_2(r),$ 
$$\xymatrix{ H_r( f^{-1}(\theta)\ar[r]^{i_1(r)} )&\overline H_r(X_\theta)  &H_r(f^{-1}\theta )\ar[l]_{i_2(r)})}.$$ 
These linear maps  define the linear relation $R(i_1(r), i_2(r)):= R^f_{\theta} (r)$ and then, by Proposition \ref{P23} the linear relation $(R^f_{\theta} (r))_\reg$
and the linear isomorphism $T^f_\theta(r): V^f_\theta(r)\to V^f_\theta(r),$ s.t. $(R^f_{\theta} (r))_\reg= R(T^f_\theta(r)).$ 

\begin{definition} \label {D2}  The $r-$ monodromy of $f:X\to \mathbb S^1$ at $\theta\in \mathbb S^1,$  for $\theta$ a weakly regular value, is the similarity class of the relation $(R^f_\theta(r))_\reg,$ equivalently  
 the similarity class of the linear isomorphism $T^f_\theta (r).$
 One denotes these similarity classes by $[(R^f_\theta(r))_\reg]$ or $[T^f_\theta (r)] .$ 
\end{definition} 
\vskip .1in

For a map  $f:X\to \mathbb S^1$   and $K$ a compact ANR denote by 
 $$\overline f_K: X\times K\to \mathbb S^1$$ the composition of $f$ with the projection of $X\times K$  on $X.$  Note that if  $\theta$ is a weakly regular  value for $f$ then it remains a weakly regular value for $\overline f_K$ and 
 $(\overline {X\times K})^{\overline f_K}_\theta= \overline X^f_\theta\times K.$ Therefore,  in view of the Kunneth formula (expressing the homology of the product of two spaces) one has

\begin{equation}\label {E5}
\begin{aligned}
V^{\overline f_{K}}_\theta (r))=& \oplus _lV^f_\theta(r-l))\otimes H_l(K)\\
T^{\overline f_K}_\theta (r) = &\oplus_l T^{f}_\theta (r-l) \otimes Id_{H_l(K)}
\end{aligned}
\end{equation}
where $ Id_{H_l(K)}$ denotes the identity map on $H_l(K).$

In particular if $K$ is contractible then  
\begin{equation}\label {E5'}
[T^{\overline f_K}_\theta (r)] = [T^{f}_\theta (r)]
\end{equation}
and if $K=\mathbb S^1$ then 
\begin{equation}\label{E6}
[T^{\overline f_K}_\theta(r)]= \begin{cases} [T^f_\theta(0)]\  \mbox {if}\  r=0\\
[T^f_\theta(r) \oplus T^f_\theta (r-1)] \  \mbox {if}\  r\geq 1.\end{cases}
\end{equation}

 \begin{proposition}\label{T31}\
\begin{enumerate}
\item  If $\theta_1$ and $\theta_2$ are two different weakly regular values for $f$ then $[T^f_{\theta_1} (r)]= [T^f_{\theta_2}(r)].$ 

\item   If $X$ is a good ANR and $f, g: X\to \mathbb S^1$ are two homotopic maps  with  $\theta_1$ a weakly regular value for $f$ and $\theta_2$  a weakly regular value for $g$ then 
 $[T^f_{\theta_1} (r)] = [T^g_{\theta_2}(r)].$ 

\item 
If $f:X\to \mathbb S^1$ and $g: Y\to \mathbb S^1$ are two maps with $\theta_1$ weakly regular value for $f$ and $\theta_2$ weakly regular value for $g$ then $[T^f_{\theta_1} (r)] = [T^g_{\theta_2}(r)]$ \ iff \ $[T^{\overline f_{\mathbb S^1}}_{\theta_1} (r)] = [T^{\overline g_{\mathbb S^1}}_{\theta_2}(r)].$

\item If $f:X\to \mathbb S^1$ and $g: Y\to \mathbb S^1$ are two maps with $\theta_1$ weakly regular value for $f$ and $\theta_2$ weakly regular value for $g,$ and $\omega:X\to Y$ is a homeomorphisms such that $g\cdot \omega$ and $f$ are homotopic then 
$[T^f_{\theta_1} (r)] = [T^g_{\theta_2}(r)].$ 
\end{enumerate}
\end{proposition}

{\it Proof of 1.: }
For $X$  a compact ANR  and  $\xi\in H^1(X;\mathbb Z)$ consider  $\pi: \tilde X\to X$ an infinite cyclic cover 
\footnote { An infinite cyclic cover is a map $\pi: \tilde X\to X$  together with  a free action $\mu:\mathbb Z\times \tilde X\to \tilde X$ such that  $\pi(\mu(n,x))= \pi(x)$ and the map induced by $\pi$ from $\tilde X/ \mathbb Z$ to $X$ is a homeomorphism. The above cover is called {\it associated to $\xi$} if any $\tilde f:\tilde X\to \mathbb R$ which satisfies $\tilde f(\mu(n,x))= \tilde f(x)+2\pi n$ induces a map from $X$ to $\mathbb R/2\pi\mathbb Z= \mathbb S^1$ representing the cohomology class $\xi_f= \xi.$ For two infinite cyclic covers $\pi_i: \tilde X_i\to X$ representing $\xi$  there exists  homeomorphisms $\omega:\tilde X_1\to \tilde X_2$ which intertwine the free actions $\mu_1$ and $\mu_2$ and  satisfy $\pi_2\cdot \omega= \pi_1.$} associated to $\xi.$

Any  map $f:X\to \mathbb S^1$ such that $f^\ast(u)=\xi,$ $u$ the canonical generator of $H^1(\mathbb S^1),$ has lifts 
$ \tilde f: \tilde X\to \mathbb R,$ which make the  diagram below with $p(t)$ is given by $p(t)= e^{it}\in \mathbb S^1,$a pull-back diagram
  
\begin{equation}
\xymatrix{ &\mathbb R\ar[r]^p &\mathbb S^1\\
&\tilde X \ar[r]^{\pi}\ar[u]_{\tilde f} &X\ar[u]_f.}
\end{equation}

 Consider $\theta_1= e^{it_1}, \theta_2= e^{it_2} \in \mathbb S^1$ two different weakly regular values for $f$ (i.e. with $t_2- t_1 \leq \pi$ hence $t_1 <t_2 <t_1+2\pi <t_2+2\pi ).$
We apply the discussion in the  subsection \ref{SS41} to the real valued map   $\tilde f: \tilde X \to \mathbb R$  and note that  

$$R^f_{\theta_1}(r) = R^{\tilde f,\tilde f}_{t_1, t_1+2\pi}(r)= R^{\tilde f, \tilde f}_{t_2, t_1+2\pi}(r)\cdot R^{ \tilde f, \tilde f}_{t_1, t_2}(r)$$ and 
$$R^f_{\theta_2}(r) = R^{\tilde f,\tilde f}_{t_2, t_2+2\pi}(r)= R^{\tilde f \tilde f}_{t_1+2\pi, t_2+2\pi}(r)\cdot R^{\tilde f, \tilde f}_{t_2, t_1+2\pi}(r).$$  

Using  the linear isomorphisms induced by  $\pi,$ the linear relations $R^{\tilde f,\tilde f}_{t_1, t_2}(r)$ and $R^{\tilde f,\tilde f}_{t_1+2\pi, t_2 +2\pi}(r)$ can be identified to the linear relation $R':= R^f_{\theta_1}(r) \colon H_r(f^{-1}(\theta_1))\rightsquigarrow H_r(f^{-1}(\theta_2)$ while $R^{\tilde f, \tilde f}_{t_2, t_1+2\pi}(r)$ to  the linear relation 
$R''= R^f_{\theta_2}(r) \colon H_r(f^{-1}(\theta_2)) \rightsquigarrow H_r(f^{-1}(\theta_2)).$ 

Therefore  $R^f_{\theta_1}(r)= R'' \cdot R'$ and $R^f_{\theta_2}(r)= R'\cdot R'',$  which in view of Proposition \ref {P23}  (5)   imply that $(R^f_{\theta_1}(r))_\reg \sim (R^f_{\theta_2}(r))_\reg.$ 
\vskip .1in 

{\it Proof of 2.:}  
In view of Proposition \ref{O52} it suffices to prove the statement under the following additional hypotheses:
\begin{enumerate}
\item  At least one  of the maps $f$ or $g$ is weakly tame,
\item $D(f,g)<\pi,$
\item The angles $\theta_1= e^{it_1}$  and $\theta_2= e^{it_2}$ satisfy $|t_2-t_1 -\pi| <\pi/4.$
\end{enumerate}
Since$f$ and $g$ are homotopic 
 $\xi_f=\xi_g.$  For any infinite cyclic cover $\tilde X \to X$ associated with $\xi= \xi_f= \xi_g$ both $f$ and $g$ have lifts 
$\tilde f$ and $\tilde g$ as indicated in the diagrams below

\begin{equation}
\xymatrix{ &\mathbb R\ar[r]^p &\mathbb S^1\\
&\tilde X \ar[r]^{\pi}\ar[u]_{\tilde f} &X\ar[u]_f} \xymatrix{ &\mathbb R\ar[r]^p &\mathbb S^1\\
&\tilde X \ar[r]^{\pi}\ar[u]_{\tilde g} &X\ar[u]_g .}    
\end{equation} 
\vskip .1in 

Under the additional hypotheses one can find lifts $\tilde f$ and $\tilde g$ such that  
$g^{-1}(t_2)\subset \tilde f^{-1}(t_1, t_1+2\pi)$    
and $\tilde f^{-1}(t_1+2\pi)\subset \tilde g^{-1}(t_2, t_2+2\pi).$  We apply the considerations  in  subsection \ref{SS41} to the real valued maps $\tilde f, \tilde g:\tilde X\to \mathbb R$  and conclude that :
$$R^{f}_{\theta_1}(r)= R^{\tilde f,\tilde f}_{t_1, t_1+2\pi}(r)= R^{\tilde g,\tilde f}_{t_2, t_1+2\pi}(r)\cdot R^{\tilde f,\tilde g}_{t_1, t_2} (r)$$ 
and $$R^{g}_{\theta_2}(r)= R^{\tilde g,\tilde g}_{t_2, t_2+2\pi}(r)= R^{\tilde f,\tilde g}_{t_1+2\pi, t_2+2\pi}(r) \cdot R^{\tilde g,\tilde f}_{t_2, t_1+2\pi}(r).$$

Let $R' := R^{\tilde g,\tilde f}_{t_2, t_1+2\pi}(r) $ and 
$R'' := R^{\tilde f,\tilde g}_{t_1, t_2}(r) = R^{\tilde f,\tilde g}_{t_1+2\pi, t_2+2\pi}(r).$
Then  $R^{f}_{\theta_1}(r)= R'' \cdot R' $ and $R^{ g}_{\theta_2}(r)= R'\cdot R''$
which,  by Proposition  \ref{P23} item (5),  imply that $(R^f_{\theta_1}(r))_\reg \sim (R^g_{\theta_2}(r))_\reg.$ 

\vskip .1in
{\it Proof of 3.:}
Recall that for a linear isomorphism $T:V\to V$ one denotes by $\mathcal J(T)$ the set of Jordan cells which is a similarity invariant.

First observe that if $T_1:V_1\to V_1$ and $T_2:V_2\to V_2$ are two linear  isomorphism then $\mathcal J(T_1\oplus T_2)= \mathcal J(T_1)\sqcup \mathcal J(T_2). $ 

If so 
 $[T_1\oplus T_2]= [T_1'\oplus T_2'],$   hence $\mathcal J([T_1]) \sqcup J([T_2]) =  \mathcal J ([T_1'])\sqcup \mathcal J([T_2']),$ 
and 
$[T_1]= [T_1'],$ hence $\mathcal J([T_1]) = \mathcal J([T_1']) ,$  imply  
$ \mathcal J([T_2])= \mathcal J([T_2']),$ hence $[T_2]= [T_2']$. 

We apply this observation to $T_1= T^f_{\theta_1}(r-1),$ $T_1'= T^g_{\theta_2}(r-1)$ and
$T_2= T^f_{\theta_1}(r),$ $T_2'= T^g_{\theta_2}(r).$
Then  by induction on $r$ formula (\ref{E6}) implies item 3.. 

\vskip .1in 
{\it Proof of 4.:} In view of item 2. one has $[T^{g\cdot \omega}_{\theta_2}(r)]= [T^f_{\theta_1}(r)].$
Since $\omega$ induces  a homeomorphism  between  $\overline X_{\theta_2}^{g\cdot \omega}$ and $\overline Y_{\theta_2}^{g \cdot \omega}$ then $R^{g\cdot \omega}_{\theta_2}(r) \sim R^{g}_{\theta_2}(r) $ which implies $[T^{g\cdot \omega}_{\theta_2}] = [ T^{g}_{\theta_2}] $ which implies  $[T^f_{\theta_1}(r)]= [T^g_{\theta_2}(r)].$

q.e.d.
\vskip .2in

In view of Proposition \ref{T31} (1) $[T^f_\theta(r)]$ is independent on  $\theta,$ so for a {\it weakly tame map} $f$ one can write 
$[T^f(r)]$ instead of $[T^f_\theta(r)].$ In view of Proposition \ref{T31} (2) if  $f_1$ and $f_2$ are two good maps with $D(f_1, f_2))<\pi$
then one has $[T^{f_1}(r)]= [T^{f_2}(r)].$ 

If $X$ is a {\it good ANR} for a map $f$ choose a weakly tame maps $f'$ with $D(f,f')<\pi/2$ and in view of Proposition \ref{T31} (2) 
 $[T^{f'}(r)]$ provides an unambiguous definition of the $r-$monodromy for the map $f.$  Indeed for two such maps $f'_1$ and $f'_2$ 
 one has $D(f'_1, f'_2))<\pi$ and then Proposition \ref{T31} (2) guaranties that  $[T^{f'_1}(r)]= [T^{f'_2}(r)].$  
 Moreover,
  if $f$ and $g$ are homotopic then $[T^f(r)]= [T^g(r)].$
Then for $X$ a {\it good ANR} and $\xi \in H^1(X;\mathbb Z)$ one chooses $f,$ with $\xi_f=\xi,$ and  one defines  $$[T^{(X;\xi)}(r)]:= [T^f(r)].$$ 

In  order to show that $[T^{(X,\xi)}(r)]$ can be extended to any compact ANR  and is a homotopy invariant of the pair $(X,\xi),$
\footnote {This means that for $(X_1,\xi_1),$ and $(X_2,\xi_2)$  pairs with $X_i, i=1,2$ compact ANRs, $\xi_i\in H^1(X_i; \mathbb Z),  i=1,2,$ the existence of a homotopy equivalence  $\omega:X_1\to X_2$  satisfying  $\omega^\ast (\xi_2)= \xi_1$ implies $[T^{(X_1,\xi_1)}]= [T^{(X_2,\xi_2)}].$}
one  uses  Proposition \ref{T31} (3) and (4) and the Stabilization  Theorem below. This theorem is  a consequence remarkable topological results of Edwards and Chapman about Hilbert cube manifolds, cf \cite{CH}. An homological proof is also possible but requires  a little bit of algebraic topology, cf \cite {BH}.
\vskip .2in 

\begin{theorem}\label{T35}{\bf Stabilization theorem} (R. Edwards and T. Chapman) 
There exists a contractible compact ANR, $Q,$ with the following  properties.

1. For any compact ANR $X$ the product 
 $X\times Q$ is a {\it good  compact ANR}. 

2. Given $\omega : X\to Y$ a homotopy equivalence of compact ANR's  the map 
\ \ \quad \quad
$\omega\times Id_{Q\times \mathbb S^1} : X \times Q \times \mathbb S^1 \to  Y \times Q \times \mathbb S^1$ is homotopic to a homeomorphism 
$\omega' : X \times Q \times \mathbb S^1 \to  Y \times Q \times \mathbb S^1.$  

The compact ANR,  $Q,$ is the product of countable many copies of the segment $[0,1].$
\end{theorem}

The statements above  are  rather straightforward  consequences of  Edwards and Chapman results  however 
neither 1. nor  2., as formulated above,  can be found in their work or in \cite {CH}. They can be derived from the  mathematics presented in \cite {CH} as explained in  Appendix 2. 

\vskip .2in 
{\bf Extension of $r-$monodromy  to all pairs $(X,\xi)$}

To  any pair $(X, \xi),$ $X$ compact ANR and  $\xi \in H^1(X;\mathbb Z),$ and any $r\in \mathbb Z_{\geq 0} $ one defines the 
$r-$monodromy  by 
$$ [T^{X,\xi}(r)] := [T^{X\times K, \overline \xi}(r)] $$  with $\overline \xi $ is the pull back of $\xi$ by the projection of $X\times K\to X.$
In view of the equality (\ref {E6}) if $X$ was already a good ANR then $[T^{X,\xi}(r)]= [T^{X\times K,\overline \xi}(r)].$  

To verify the homotopy invariance consider $f_i:X_i\to \mathbb S^1$ representing the cohomology class $\xi_i.$  Since $\omega^\ast (\xi_2)= \xi_1$ the composition $f_2\cdot \omega\ $ and $f_1$ are homotopic and then in view of item 2. of Stabilization Theorem one has the homeomorphism $\omega'$  homotopic to $\omega\times id_{K\times \mathbb S^1}.$ Hence $(\overline {f_2})_{K\times \mathbb S^1} \cdot \omega'$ is homotopic to $(\overline {f_1})_{K\times \mathbb S^1}.$ Then by  Proposition \ref{T31} (4) one has
$[T^{(\overline {f_2})_{K\times \mathbb S^1}}(r)]= [T^{(\overline {f_1})_{K\times \mathbb S^1}}(r)].$  In view of Proposition  \ref{T31} (3), \  $[T^{(\overline {f_2})_{K}}(r)]= [T^{(\overline {f_1})_{K}}(r)]$, hence 
$[T^{(X_1, \xi_1)}]= [T^{(X_2, \xi_2)}].$

As a summary  one has the following Theorem.   

\begin{theorem} \label {T36} To any pair $(X,\xi)$, and $r =0,1,2, \cdots,$  $X$ compact ANR,  and $\xi \in H^1(X;\mathbb Z)$ one can associate the similarity class of linear isomorphisms $[T^{(X,\xi)} (r)]$ which is a homotopy invariant of the pair.
When $f:X\to \mathbb S^1$  is a  good map with $\xi_f=\xi$  this is the $r-$monodromy  defined for a good map $f$  and a weakly regular value.  
\end{theorem} 

The collection $\mathcal J_r(X;\xi)$ consisting of the pairs  with multiplicity, 
$(\lambda, k), \ \lambda \in \overline \kappa, k\in \mathbb Z_{>0},$  which determine  the similarity class $[T^{(X;\xi)}(r)]=  
[\oplus _{(\lambda,k)\in \mathcal J_r(\xi)} T(\lambda, k)]$  is referred to as {\it the Jordan cells of the $r-$ monodromy} of $(X,\xi).$ 
\vskip .2in 

{\bf An example}

\noindent The picture below  is taken from \cite{BD11} but with a different gluing map.

\begin{figure}[h]
\includegraphics{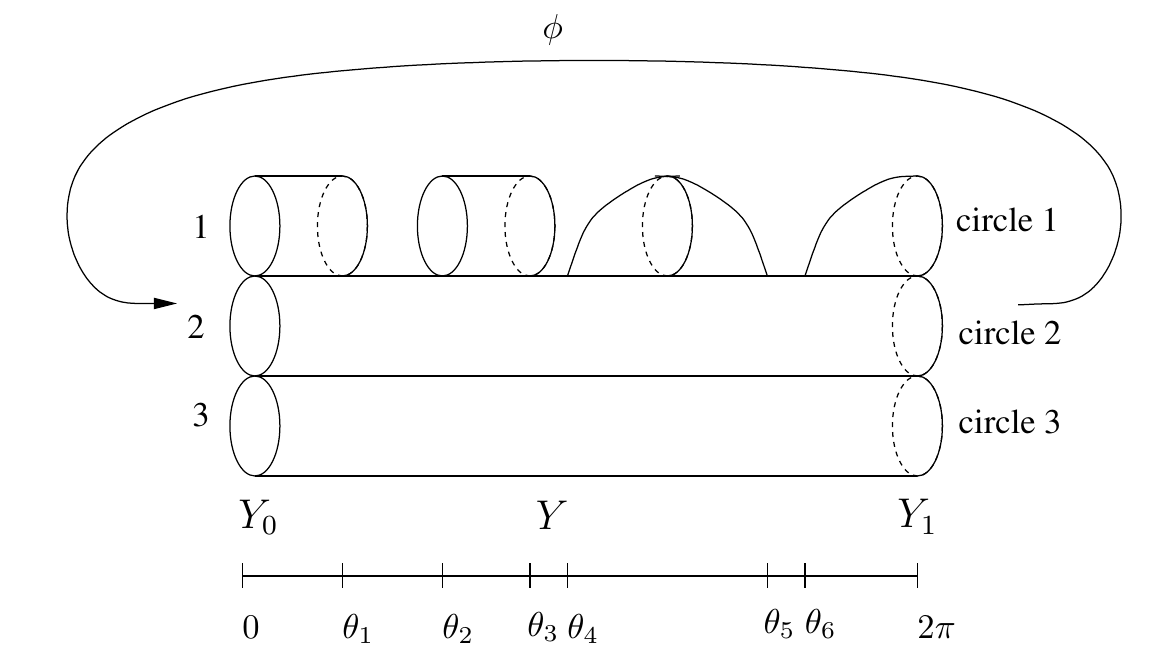}
\end{figure}

Consider the  space $X$ obtained from $Y$  by identifying its  
right end $Y_1$ (a union of three circles) to the left end 
$Y_0$ (a union of three circles) following 
the map $\phi\colon Y_1\to Y_0$ defined by the matrix 
\begin{equation*}
\begin{pmatrix}
3&3 & 0\\
2&3&-1\\
1 & 2&3            
\end{pmatrix}.
\end{equation*}
The meaning of this matrix as a map is the following: Circle (1)  is divided in 6 parts, circle (2) in 8 parts  and and circle (3) in 4 parts  ; the first three parts of circle (1) wrap  clockwise  around  circle (1) to cover it  three times,  the next 2 wrap clockwise around circle (2) to cover it twice and around circle three to cover it three times.  Similarly circle (2) and (3) wrap over circles (1)(2) and (3) as indicated by the matrix. The first part of circle (2) wraps counterclockwise on circle (2).

The map $f\colon X\to S^1$ is induced by 
the projection of $Y,$ on the interval $[0,2\pi]$ which becomes $\mathbb S^1$ when $0$ and $2\pi$ are identified. This map has all values weakly regular. 

In this example  $\mathcal J_0(f)= \{ (\lambda=1,k=1)\},$ $\mathcal J_1(f)= \{(\lambda=2, k=2)\}$ and $\mathcal J_2(f)= \emptyset.$  The first and last calculations are obvious. The second   will be derived by applying the algorithm described in Section \ref{S4}.

\subsection {\bf F- monodromy }

For a field $\kappa,$ instead of the homology vector space $H_r(X),$  one can consider a more general functor $F,$ a so called Dold half-exact functor cf \cite {D}.
Recall that this is a covariant functor defined from  the category $Top_c$ of compact ANR's and continuous maps (or any subcategory with the same homotopy category) to the category $\kappa-Vect$ of finite dimensional vector spaces and linear maps which satisfies the following properties:
\begin{enumerate}
\item $F(f)=F(g)$ for any two homotopic maps $f$ and $g,$
\item $F$ satisfies the Meyer Vietoris property.  Precisely, 
if $A$ is a compact ANR with $A_1$ and $A_2$ closed subsets such that $A_1, A_2$ and $A_{12}= A_1\cap A_2$ are ANR's and $A=A_1\cup A_2$ then the sequence 
$$\xymatrix{ F(A_{12})\ar[r] ^i & F(A_1)\oplus F(A_2)\ar[r]^j & F(A)}$$ is exact. Here   $i= F(i_1) \oplus F(i_2),$   $j= F(j_1) - F(j_2),$
$i_1, i_2$ the obvious inclusions of $A_{12}$ in $A_1$ and $A_2$ and $j_1, j_2$ the obvious inclusion of $A_1$ and $A_2$ in $A.$
\end{enumerate} 
Analogues of Propositions \ref{T30} and  \ref {T31} hold for $F$ instead of $H_r$ since they are based only on the Meyer-Vietoris property. 

The same constructions applied to $F$ instead of $H_r$ work and one defines the $F-monodromy$ on the same lines.
There are plenty of such functors and the $F-$monodromy might be an  invariant which deserve attention.

\subsection{\bf  Comments }\label{ss34}

Theorem \ref{T36} is implicit  in \cite{BH} (cf section 4 combined with  with Theorem 8 .14)  based on the interpretation of the monodromy as the similarity class of the linear isomorphism induced by the generator of the group of deck transformations,
on the vector space $V_r(X,\xi)= \ker (H_r (\tilde X)\to H^N_r(X,\xi)).$ 
Here  $\tilde X$ denotes is the infinite cyclic cover of $X$ defined by $\xi$ and $H^N_r(X;\xi)$ denotes the Novikov homology of $(X,\xi).$ 
\vskip .1in

In \cite{BH} it is shown that the Jordan cells $\mathcal J_r(f)$ defined in \cite{BD11} as invariants for persistence of the circle valued map $f$  are the same as the Jordan cell defined above. 
Since \cite{BH} is not yet in print, for the reader familiar with the notations in  \cite {BD11} section 5, we will provide a short explanations of this statement in Appendix 3. 
\vskip .1in

The characteristic polynomial of $[T^{(X,\xi)}(1)] $ for the pair $(X;\xi),$  $X=S^3\setminus K,$ $K$ an open tube around an embedded oriented circle (knot) and $\xi$ the canonical generator of $H^1(S^3\setminus K)= \mathbb Z$ is exactly the  monic Alexander polynomial of the knot.
as explained in section \ref{S5}.
\vskip .1in

The alternating product of the characteristic polynomials $P_r(z)$ of the monodromies $[T^{X;\xi}(r)],$
$$A(X; \xi)(z) = \prod P_{r}(z)^{(-1)^r},$$   known to topologists  as the Alexander rational function, calculates (essentially \footnote {a precise formulation requires additional data whose explanations are beyond the purpose of this paper}) the Reidemeister torsion of $X$ equipped with the degree one representation  of $\pi_1(X)$ defined by $\xi,$ and the complex number $z\in \mathbb C, \ z\ne 0$ 
when interpreted as an homomorphism $\pi_1(X,x)\to GL_1(C).$  This was pointed out first by J. Milnor and refined by V. Turaev, cf \cite{Tu}. 
Given the need  of additional background and definitions  a precise formulation of this calculation  will be discussed elsewhere.  

\section {The  calculation of Jordan cells of an angle valued map}\label{S4}

\subsection {\bf Generalities}

{\bf Cell complexes}

Recall that:
\begin{itemize}
\item
A {\it convex $k-cell$}  $\sigma$ in an affine space $\mathbb R^n, n\geq k,$ is the convex hull of a finite collection of points  $e_0, e_1, \cdots e_N$ called vertices,   with the property that :
\begin{enumerate}
\item 
there are subsets  with $(k+1)-$points  linearly independent \footnote { In an affine space $(k+1)-$ points are linearly independent if they lie in a $k-$dimensional affine subspace but not in any $(k-1)-$dimensional affine subspace.} 
but no subset of $(k+2)-$points linearly independent,   
\item 
no vertex  lies in the topological interior of this convex hull 
\end{enumerate}
The topology of the cell is the one induced from the ambient affine space $R^n.$

A $k-$ simplex is a convex $k-$ cell with exactly $k+1$ vertices.

A  $k'-$face  $\sigma'$ of $\sigma,$ $k'<k,$ is a convex $k'-$cell  whose set of vertices is a subset of the set of vertices of $\sigma.$ One  indicates that $\sigma'$ is a face of $\sigma$ by  writing $\sigma'\prec \sigma.$ 

A space homeomorphic to a convex $k-cell$ is called simply a $k-cell$ and the subset homeomorphic to a face continues to be called {\it face}.  

\item A  {\it cell  complex} $Y$ is a space together with a locally finite  collection $\mathcal Y$  of compact subsets $\sigma\subset Y,$ each  homeomorphic with a convex cell  with the following properties:
\begin{enumerate}
\item  If a $k-$ cell $\sigma$ is a member of the collection $\mathcal Y$ then any of its faces $\omega \prec \sigma$ is a member of the collection $\mathcal Y.$
\item If $\sigma$ and $\sigma'$ are two cells members of the collection $\mathcal Y$ then their intersection is a union of cells and each cell of this union is  face of both $\sigma$ and $\sigma'.$
\item $\bigcup_{\sigma\in \mathcal Y} \sigma= Y.$
\end{enumerate}

The concept of sub complex $Y'\subset Y$ is obvious. Precisely $Y'$ is the union of the cells in the  subset  $\mathcal Y' \subset \mathcal Y$ with the property that any face of cell in $\mathcal Y'$ is in $\mathcal Y'.$ 

A finite {\it simplicial complex} is a finite cell complex with all cells simplexes.

For a cell complex $Y$ with cells $\mathcal Y$ denote by $\mathcal Y_k$ the set of the $k-$ cells in $\mathcal Y.$ Clearly $\mathcal Y_0$ is the set of all vertices of the cells in $\mathcal Y.$  

\item An {\it oriented cell} is a cell $\sigma\in \mathcal Y$ equipped with an orientation  $o(\sigma).$ This orientation induces an orientation for any codimension one face $\sigma$ described by the rule : {\it first the induced orientation, next the normal vector pointing inside give  the orientation $o(\sigma)$}.  

If each cell $\sigma$ of a cell complex is equipped with an orientation $o(\sigma)$
one has the incidence function $\mathbb I:\mathcal Y\times \mathcal Y  \to \{0,+1, -1\}$ defined as follows: 

\begin{equation}\label{E100}
\mathbb I(\sigma, \tau):= 
\begin{cases}
\mathbb I(\tau, \sigma)= +1 &\ if \  \sigma\in \mathcal Y_k, \tau \in \mathcal Y_{k+1}, \sigma \prec \tau, o(\sigma)|_{\sigma'}= o(\sigma'),\\
\mathbb I(\tau, \sigma)=  -1 &\ if \  \sigma\in \mathcal Y_k, \tau \in \mathcal Y_{k+1}, \sigma \prec \tau, o(\sigma)|_{\sigma'}\ne o(\sigma'),\\ 
\mathbb I(\tau, \sigma)=   0  &\ if \  \sigma\cap \sigma'=\emptyset
\end{cases}
\end{equation}

The incidence function  determines the homology of $Y$ with coefficients in any field. 

\item A {\it good order} of the set $\mathcal Y$ of cells of $Y$  is a total order ''$ \leq $ " if: 

(1) $\sigma \prec \tau $ implies $\sigma <\tau .$

In this case, if the cardinality of $\mathcal Y$ is $N,$ then the function $\mathbb I(\cdots, \cdots)$ can be regarded as  $N\times N$ upper triangular matrix (all entries on and below diagonal are $0$ ) and  is referred below as the {\it incidence matrix} of $Y.$ 

Suppose  that inside $Y$ one has two disjoint sub complexes, $Y_1, Y_2\subset Y.$
In this case 
a  {\it good order 
compatible with $\mathcal Y_1$ and $\mathcal Y_2$}  needs,
in addition to (1) above, the following requirements  satisfied:

(2) If $\sigma_1\in \mathcal Y_1$ and $\sigma_2 \in \mathcal Y_2$ then 
$\sigma_1 \prec \sigma_2$  and

(3) If $\sigma' \in \mathcal Y_i$ and $\sigma \in \mathcal Y\setminus \mathcal Y_1\sqcup \mathcal Y_2$ imply $\sigma' \prec \sigma.$ 

\end{itemize}

Note that: 
\begin{enumerate}
\item Given a  total order of the cells in $\mathcal Y$ a simple algorithm referred below as the {\it Ordering algorithm} permits to  correct  this order into a good total order. 
The  Ordering Algorithm  is based on the inspection of the $n-$th cell  with respect to all previous cells. If  the requirements (1)-(3) are not violated move to the $(n+1)-$cell. If at least one of the three requirements is  violated, change the position  of this cell, and implicitly of the preceding ones if the case, by moving  the cell to the left until (1), (2), or (3) are no more violated.

\item With the requirements 1, 2, 3  for the  {\it good order} satisfied  the incidence matrix of $Y,$ $\mathbb I(\cdots, \cdots),$ has the form 

\begin{equation}
\begin{pmatrix}
A_1&  0&  X\\
0& A_2&  Y\\
0&  0&  Z
\end{pmatrix}
\end{equation}
with $A_1= \mathbb I_1$, $A_2= \mathbb I_2$ the incidence matrices for  $Y_1$ and for $Y_2.$

\item
The persistence algorithm \cite{CEM}, \cite{ZC} permits to calculate from the incidence matrix $\mathbb I$
\begin {enumerate}
\item first, a base for $H_r(Y_1),$ then a base for $H_r(Y_2),$ then a base for $H_r(Y),$ 
\item second,  the $\dim H_r(Y) \times \dim H_r(Y_1)-$matrix $A$ and the 
\newline $\dim H_r(Y) \times \dim H_r(Y_2)-$matrix $B$ 
representing the linear maps induced in the  homology in dimension $r$ by the inclusions of $Y_1$ and $Y_2$ in $Y.$
\end{enumerate}
\vskip .1in
\end{enumerate}

{\bf The cut at $\theta$  
(with respect to $f$).}  

Let  $\sigma$ be  a $k-$dimensional simplex with vertices $e_0, e_1, \cdots e_k,$  i.e. a convex $k-$cell generated by $(k+1)$ linearly independent points located in some vector space.  Let  $f:\sigma\to \mathbb R$ be a linear map determined by the values of $f(e_i)$  by the formula 
\begin{equation}\label{E12}
 f(\sum_{i} t_i e_i)= \sum_{i} t_i f(e_i), \  t_i\geq 0 , \sum_i t_i=1 
\end{equation}
and let $t\in \mathbb R.$  
Suppose that  $\sup_i f(e_i) >t$ and $\inf_i f(e_i) <t.$ 

The map  $f$ and the number  $t$ determine two  $k-$convex cells $\sigma_+, \sigma_-$ and a $(k-1)-$convex cell $\sigma'$:
\begin{equation}\label {E18}
\begin{aligned}
\sigma_+= &f^{-1}([t, \infty))\cap \sigma \\
\sigma_-= &f^{-1}((-\infty, t])\cap \sigma \\
\sigma'= &f^{-1}(t)\cap \sigma  .
\end{aligned}
\end{equation}

An orientation  $o(\sigma)$ on $\sigma$ provides orientations $o(\sigma_+), o(\sigma_-)$on $\sigma_+,\sigma_-$ and induces an orientation $o'(\sigma')$ on $\sigma',$ precisely the unique orientation which followed by  the direction provided by the vector field $\mbox {grad} f$ defines  the orientation of $o(\sigma).$  Then $I(\sigma_\pm, \sigma')= \pm 1.$

Recall that the map $f: X\to \mathbb S^1\subset \mathbb C$ is simplicial if the restriction of ''$ \ln f$''  \footnote {in view of 1-connectivity of each simplex  "$\ln f$"  has continuous univalent determination when the value on one vertex of the simplex is specified} to any simplex $\sigma$ is a linear map as considered above.

\subsection {\bf The algorithm} \label {SS42}

The algorithm we propose inputs a simplicial complex  $X,$ a simplicial map $f$ and an angle $\theta$  different from the values of $f$ on vertices and outputs for each $r,$ $0\leq r\leq \dim X,$ in STEP 1 
two  $m\times n$ matrices $A_r$ and $B_r$ with $m,$ the number of rows equal to the dimension of $H_r(\overline X^f_\theta)$ and $n,$ the number of columns, equal to the dimension of $H_r(f^{-1}(\theta)).$ 
The matrices represent the  linear maps induced in homology by the two inclusions, from the left and from the right,  of $f^{-1}(\theta)= Y_1 = Y_2$ in $\overline X^f_\theta.$
In STEP 2 one obtains from the matrices $A_r$ and $B_r$ the invertible square matrices $A'_r$ and $B'_r$ such that  $(B'_r)^{-1} A'_r$ represents the $r-$monodromy and in STEP 3 one derives from $(B'_r)^{-1} A'_r$ the Jordan cells $\mathcal J_r(f).$ 
\vskip .1in
{\bf STEP 1.}

The simplicial set  $X$ is recorded by :

-- the set of  vertices with an arbitrary  chosen total order, 

-- a  specification of the subsets which define the collection $\mathcal X$ of simplices.   

Implicit in this data is an orientation $o(\sigma)$ of each simplex, orientation provided by the relative ordering of the vertices of each simplex,  and therefore the incidence number $\mathbb I(\sigma', \sigma)$  of any two simplexes $\sigma'$ and $\sigma$ in $\mathcal X.$

(Implicit is also a total order of the simplexes of $\mathcal X$ provided by the {\it lexicographic order}  induced from the order of  the vertices.)

The simplicial map $f$ \footnote {for simplicity one supposes that $f$ takes different values on different vertices} is indicated by 

-- the sequence of $N_0=$ the number of vertices, the values of $f$ on vertices. 

The map $f$ and the angle $\theta=  e^{it}$ provide a decomposition of the set $\mathcal X$ as  $\mathcal X'\sqcup \mathcal X''$ with $\mathcal X':= \{\sigma\in \mathcal X \mid \sigma \cap f^{-1}(\theta)\ne \emptyset  \} $ and $\mathcal X'':=\mathcal X\setminus \mathcal X".$

From these data we can derive :

-- first,  the collections $\mathcal Y$ with the  sub collections $\mathcal Y(1)$ and $\mathcal Y(2)$  of the cells of the complex $Y= \overline X^f_\theta$ and the sub complexes $Y_1= f^{-1}(\theta)$ and $Y_2= f^{-1}(\theta),$ 

-- second, the  incidence function on $\mathcal Y\times \mathcal Y,$ 

-- third,  a  good order for the elements of $\mathcal Y.$ 

These all lead to the incidence matrix $\mathbb I (Y).$ 

{\it Description of the cells of $Y:$}  Each oriented simplex $\sigma$ in $\mathcal X''$ provides a unique oriented cell $\sigma$  
in $\mathcal Y.$ 

Each oriented $k-$simplex $\sigma$ in $\mathcal X'$ provides two oriented $k-$cells $\sigma_+$ and $\sigma_-$ and two oriented 

$(k-1)-$cells $\sigma'(1)$ and $\sigma'(2),$ copies of the oriented cell $\sigma'.$ So the cells of $Y$ are of five types 

$\mathcal Y'_k(1)= \mathcal X'_{k+1},$ 

$\mathcal Y'_k(2)= \mathcal X'_{k+1},$ 

${\mathcal Y'_k}_-= \mathcal X'_{k}, $

${\mathcal Y'_k}_+= \mathcal X'_{k},$  

$\mathcal Y_k''= \mathcal X''_{k}.$

Note that ${\mathcal Y'_k}_+ $ and 
${\mathcal Y'_k}_-$ are two copies of the same set $ \mathcal X'_k$  and  $\mathcal Y'_k(1)$ and 
$\mathcal Y'_k(2)$ are in bijective correspondence  with the set  $\mathcal X'_{k+1}.$

Inside the cell complex $Y$  we have two sub complexes $Y_1$ and $Y_2$  whose cells are 
$ (\mathcal Y_1)_k = \mathcal Y'_k(1), $
$ (\mathcal Y_2)_k = \mathcal Y'_k(2), $  two copies of the same set $\mathcal X'_{k+1}.$
\vskip .1in
{\it Incidence of cells of $\mathcal Y:$} 
The incidence of two cells in the same group (one of the five types) are the same as the incidence of the corresponding simplexes. 
The incidence of  two cells one in  $\mathcal Y_1$ the other in  $\mathcal Y_2$ or  one in the group $Y'(i), i=1,2$  the other in  the group $\mathcal Y''$ is always zero. The  rest of  incidences are provided by the formulae (\ref{E100}).
\vskip .1in
{\it The good  order: } Start with a good  order of  $\mathcal Y_1$  followed by  $\mathcal Y_2$ with the same order  (translated by the number of the elements of $\mathcal Y_1$) 
followed by the remaining elements of $\mathcal Y.$  Without changing the order in  the collection $\mathcal Y_1\sqcup \mathcal Y_2,$ 
since no violation of the requirements 1, 2, 3, appear,  one can realize a good order for the entire collection $\mathcal Y$ 
with all remaining cells being preceded by the cells of $\mathcal Y_1 \cup \mathcal Y_2.$ Simply we apply the {\it ordering algorithm } to obtain a good order. 

As a result we have the incidence matrix $\mathbb I(Y)$ which is of the form 

\begin{equation}
\begin{pmatrix}
\mathbb I&  0&  X\\
0&  \mathbb I&  Y\\
0&  0&  Z
\end{pmatrix}
\end{equation}
with $\mathbb I$ the incidence matrix of $Y_1$ and $Y_2.$ 

Running the persistence algorithm cf \cite {CEM}, \cite {ZC}
leads to the matrices representing $A_r: H_r(Y_1)\to H_r(Y)$ and $B_r: H_r(Y_2)\to H_r(Y)$ as follows. 

We run the persistence algorithm  on the incidence matrix $A$ to compute a base for of the homology of $H_r(Y_1)= H_r(Y_2)$ . We continue the procedure by adding columns and rows to the matrix to obtain a base of $H_r(Y).$ It is straightforward to compute a matrix representation for the the inclusion induced linear maps $H_r(Y_i)\to H_r(Y), i=1,2.$

The time complexity of this step was discussed in \cite {BD11} and is $O(M(n)),$  the time complexity of multiplying $n\times n$ matrices, with $M(n)= O(n^\omega),$ $\omega<  2.376$  see references in \cite{BD11}.
 \vskip .1in
 
 {\bf STEP  2.} One uses the algebraic algorithm to pass from $A_r, B_r$ to the invertible matrices $A'_r, B'_r$ and then to $(B'_r)^{-1}  (A')_r$  described in the next subsection. This is based by reducing matrices to echelon form as described in subsection \ref {SS42} below.
\vskip .1in

{\bf STEP 3.} One uses the standard algorithms to put the matrix $(B')_r^{-1}  A'_r$ in Jordan diagonal form (i.e. as block diagonal matrix with Jordan cells on diagonal). 
Since the resulting matrix from Step 2 is a $k\times k$ matrix with $k<\inf\{m,n\}$ the time complexity of STEP 3  is  at most $O( k^5 \log k),$ 
cf \cite {PM}, however there are apparently  algorithms with better time complexity like \cite {Chen} with $O(5/4 k^4) $ and I understand even $O(M(k)).$
 
All basic softwares which carry linear algebra packages  contain sub packages which input a matrix and output its (reduced) row/column echelon form as well as the  matrix $C$ or $D$ in Proposition \ref{P6} involved in Step 2;  this permit an easy implementation of step 2.    
Most of them  also contain sub packages  which input a square matrices and output their Jordan form  making also easy the implementation of step 3.

\subsection {\bf An algorithm for the calculation of $R(A,B)_\reg$}\label {SS42}

 The algorithm presented below (STEP 2.) above inputs   
 two  $m\times n$ matrices $(A,  B)$ defining a linear relation $R(A,B)$ 
and outputs  two $k \times k , k\leq \inf \{m,n\},$ invertible matrices $(A', B')$
such that $R(A,B)_\reg\sim  R(A', B')_\reg.$  It is based on three modifications $T_1, T_2, T_3$ described below.
The simplest way to perform these  modification  is to use familiar procedures  of bringing a matrix to row or column echelon form 
(REF) or (CEF) explained below, but  less is actually needed as the reader can see in the presentation of the algorithm. 
\vskip .1in
{\bf Modification $T_{1} (A,B)= (A', B')$:}  

 Produces the invertible $m\times m$ matrix  $C$ 
and the invertible  $n\times n$ matrix $D$ so that 

$CAD= \begin{pmatrix} A_{11}& A_{12} \\0& 0\end{pmatrix}$ and 
$CBD= \begin{pmatrix} B_{11}& B_{12} \\B_{2.1}& 0\end{pmatrix}.$

Precisely, one constructs first $C$ which puts $A$ in REF (reduced row echelon form) 
such that  

$CA = \begin{pmatrix} A_{1}\\
0 \end{pmatrix}$ 
and makes
$CB = \begin{pmatrix} B_{1}\\
B_2 \end{pmatrix}.$

Second, one constructs $D$ which puts $B_2$ in CEF ( column echelon form). Precisely, 

$B_2  D = \begin{pmatrix} B_{21}& 0
\end{pmatrix}.$  

Clearly $CAD$ and $CBD$ are as stated above.

Take $A'= A_{12}, B'=B_{12}.$

In view of Proposition \ref{P24} (1) one has $R(A,B)_\reg= R(A',B')_\reg.$

\vskip .1in

{\bf Modification $T_{2} (A,B)= (A', B')$:}   

Produces the invertible $m\times m$ matrix  $C$ 
and the invertible $n\times n$ matrix $D$ so that 

$CAD= \begin{pmatrix} A_{11} &A_{12} \\ A_{21} , &0\end{pmatrix}$ and 
$CBD= \begin{pmatrix} B_{11} &B_{12}\\ 0  &0\end{pmatrix}.$

Precisely, one constructs $C$ which puts $B$ in REF (row echelon form) 
such that  

$CB = \begin{pmatrix} B_{1}\\
0 \end{pmatrix}$ 
and  makes
$CA = \begin{pmatrix} A_{1}\\
A_2 \end{pmatrix}.$

Then one constructs $D$ which puts $A_2$ in RCEF ( column echelon form ),
precisely 
$A_2  D = \begin{pmatrix} A_{21}& 0
\end{pmatrix}.$ 

Take $A'= A_{12}, B'=B_{12}.$

Clearly $CAD$ and $CBD$ are as stated above.

In view of Proposition \ref{P24} (1) one has $R(A,B)_\reg= R(A',B')_\reg.$

Note that if $A$ was surjective then $A'$ remains surjective.

\vskip .1in 
{\bf Modification $T_3 (A,B)= (A', B')$:} 

Produces the invertible $n\times n$ matrix $D$  and the $m\times m$ invertible matrix $C$ 
so that 

$CAD= \begin{pmatrix} A_{11}& 0 \\A_{21} & 0\end{pmatrix}$ and 
$CBD= \begin{pmatrix} B_{11}& B_{12} \\B_{21} &0\end{pmatrix}.$

Precisely, one constructs $D$ which puts $A$ in CEF (reduced row echelon form) 
i.e. 

$AD = \begin{pmatrix} A_1 &
0 \end{pmatrix}$ 
and  makes  
$BD= \begin{pmatrix} B_1&
B_2 \end{pmatrix}.$

Then one constructs $C$ to put $B_2$ in REF 
precisely,

$CB_2 = \begin{pmatrix} B_{21}\\
0\end{pmatrix}.$

Take $A'= A_{21}, B'=B_{21}.$

Clearly $CAD$ and $CBD$ are as stated above.

In view of Proposition \ref{P24} (2)  one has  $R(A,B)_\reg= R(A',B')_\reg.$

Note that if both $A$ and $B$ were surjective then  $A'$  and  $B'$  remain surjective.

\vskip .1in
Here is how the algorithm works. 

\begin{itemize}

\item  (I) Inspect $A$ 

\hskip .1in if surjective move to  (II) 

\hskip .1in else: 

\hskip .2in  - apply $T_1$  and obtain $A'$ and $B'.$ 

\hskip .2in  - make $A=A'$ and $B=B'$ and 

\hskip .2in  - go to (I)

\item (II) Inspect $B$ 

\hskip .1in if surjective move to (III) 

\hskip .1in  else :

\hskip .2in - apply $T_2$  and obtain $A'$  and $B'.$ 

\hskip .2in   - make  $A=A'$ and $B=B'$ and 

\hskip .2in   -go to (II)

(Note that if $A$ was surjective by applying $T_2,$ $A'$ remains surjective.) 

\item (III) Inspect $A$ 

\hskip .1in if injective go to (IV).

\hskip .1in  else 

\hskip .2in  -apply $T_3$  and obtain $A'$ and $B'.$ 

\hskip .2in  - make  $A=A'$ and $B=B'$ and

\hskip .2in   - go to (III)

\item (IV) Calculate $B^{-1}\cdot A.$

(Note that if $A$ and $B$ were surjective by applying  $T_3,$  $A'$ remains surjective.)  

\end{itemize}

{\bf Echelon form for $n\times m$ matrices }

Let $\kappa$ be a field. 

Let $M$ be a  $m\times n$ matrix with coefficient in the field $\kappa .$
\begin{equation*}
M=
\begin{pmatrix}
a_{1,1} & a_{1,2}       & a_{1,3}      & \cdots  & a_{1,n}     \\
a_{2,1} & a_{2,2}       & a_{2,3}      & \cdots  & a_{2,n} \\
a_{3,1} & a_{3,2}       & a_{3,3}      & \cdots  & a_{3,n} \\
\vdots  & \ddots  & \ddots & \ddots & \vdots      \\
a_{m,1} & a_{m,2}       & a_{m,3}      & \cdots  & a_{m,n}
\end{pmatrix}_.
\end{equation*}

Recall that:

 a row or column is zero-row or zero-column  if all entries are zero,
 
 the leading entry in a row or a column is the first nonzero entry w.r.to the index which varies. 

\begin{definition}\
\begin{enumerate}
\item The matrix $M$ is in {\it 
row echelon form}, REF, if  
the following hold:

\begin{enumerate}
\item All zero rows are below nonzero ones.
\item For each row the leading entry is to the right of the leading entry of the previous row.
\end{enumerate}
As example the matrix $M$ below is in 
row echelon form 

\begin{equation*}
M=
\begin{pmatrix}
0 & 0 & \alpha & x & x & x & 0 & x \\
0 & 0 & 0 & \beta & x & x & x & x \\
0 & 0 & 0 & 0 & 0 & 0 & \gamma & x \\
0 & 0 & 0 & 0 & 0 & 0 & 0 & 0 \\
0 & 0 & 0 & 0 & 0 & 0 & 0 & 0
\end{pmatrix}_.
\end{equation*}
with $\alpha, \beta, \gamma \ne 0 $ and $x$ various (possibly zero) elements in $\kappa.$

\item The matrix $M$ is in {\it 
column echelon form},  CEF, iff the transposed matrix $M^t$ is in REF i.e.   
the following hold:
\begin{enumerate}
\item All zero columns succeed nonzero ones.
\item For each column the leading entry is below of the leading entry of the previous column.
\end{enumerate}
As example the matrix below is in reduced row echelon form 
\begin{equation*}
M=
\begin{pmatrix}

0 &  0 & 0 & 0 & 0 \\
\alpha & 0 & 0 & 0 & 0 \\
x & \beta  & 0 & 0 & 0\\ 
x & x & 0 & 0 & 0 \\
x & x & 0 & 0 & 0 \\
x & x & \gamma  & 0 & 0\\
\end{pmatrix}_.
\end{equation*}
with $\alpha, \beta, \gamma \ne 0 $ and $x$ various (possibly zero)  elements in $\kappa.$
\end{enumerate}
\end{definition}

\begin{proposition}\label {P6}\ 
\begin{enumerate}
\item  For any $(m\times n)$ matrix $M$ one can produce an invertible $n\times n$ matrix $C$ such that  the composition $C M$ is in REF.
\item  For any $(m\times n)$ matrix matrix $M$ one can produce an invertible $m\times m$ matrix $D$ such that the composition $M D$ in in CEF.
\end{enumerate}
\end{proposition}

The construction of $C$ is based on "Gauss elimination" procedure  consisting in operation of "permuting rows , multiplying rows with a nonzero element in $\kappa$  and replacing a row by itself plus a multiple of an other row,  each such operation is realizable by left multiplication by  elementary matrix or permutation matrix cf \cite {G}. 

The construction of $D$ is done by : transpose,  then apply the construction of $C,$ then transpose again. 

The complexity of applying any of the modifications $T_1$, $T_2,$ or $T_3$ , construction the matrices which perform is  reduction to echelon form  is $O(k^\omega),$ for $k= \sup \{n,m\}$  cf \cite {KS}. 

\subsection {\bf An example}

We illustrate Step 2 of the algorithm with $A= A_1$ and $B=B_1$ representing the inclusion induced linear maps in homology in dimension one  derived from  the example in Section 3. We take the cut  at the angle $\theta= 0,$ i.e. the level corresponding  to the complex number $1\in \mathbb S^1.$ One has $H_1(f^{-1}(\theta)) = \kappa^3,$  $H_1(\overline X^f_{\theta})= \kappa^4.$ It is immediate from the description of $f$  that the matrices $A_1$ and $B_1$ are  given by

\begin {equation}
A_1= A= \begin{pmatrix}  3 &3 &0  \\2 &3 &-1  \\1 &2 &3  \\0 &0 &0  \end{pmatrix} \ \  B_1= B= \begin{pmatrix}  0 &0 &0& \\0 &1 & 0  \\0 &0 &1  \\ 0 &0 &0  \end{pmatrix}. 
\end{equation}

Proceed according to the algorithm:

Inspect $A,$ since not surjective apply $T_1$ and find  $C= Id$ and $D=Id.$ 
Then 

\begin {equation}
A'= \begin{pmatrix}  3 &3 &0  \\2 &3 &-1 \\1 & 2& 3  \end{pmatrix} \ \  B'= \begin{pmatrix}  0 &0 &0  \\0 &1 & 0 \\0 &0 &1  \end{pmatrix} 
\end{equation}
Update 

\begin {equation}
A= \begin{pmatrix}  3 &3 &0  \\2 &3 &-1 \\1 & 2& 3  \end{pmatrix} \ \  B= \begin{pmatrix}  0 &0 &0  \\0 &1 & 0 \\0 &0 &1  \end{pmatrix} 
\end{equation}

Since $A$ is surjective inspect $B.$ Since $B$ is not surjective apply $T_2$ and  find 

$C= \begin{pmatrix}  0 &1 &0  \\0 &0 &1 \\1 & 0& 0  \end{pmatrix}$
and then 
 $D= \begin{pmatrix} 1 &-1 &0  \\0 &1 &0 \\0 & 0& 1  \end{pmatrix}.$

Then $CAD= \begin{pmatrix} 2 &1 &-1  \\1 &1 &3 \\3 & 0& 0  \end{pmatrix}$          and $CBD= \begin{pmatrix} 0 &1 &0  \\0 &0 &1 \\0 & 0& 0  \end{pmatrix}$

\begin {equation}
A'= \begin{pmatrix}  1 &-1  \\1 &3  \end{pmatrix} \ \  B'= \begin{pmatrix}  1 &0   \\0 &1 .  \end{pmatrix} 
\end{equation}
Since both 
$A'$ and $B'$ are invertible,  consider 
\vskip .1in 
\hskip 1 in  $B^{-1}\cdot A= \begin{pmatrix}  1&-1\\1&3     \end{pmatrix}$
invertible matrix.
\vskip .1in
According to Step 3.  $\mathcal J ([R(A,B)_{\reg}])= \{(2,2)\}.$

\section {A few computational applications} \label {S5}

In view of the previous section we view the Jordan cells as {\it computer friendly} invariants (since they are computable by implementable effective algorithms). In this section we indicate how they can be used, possibly complemented by other computer friendly invariants like the standard Betti numbers 
to derive relevant topological invariants  of interest  even outside topology. 
More details and additional 
computations (applications) will be treated in future work. 

\subsection {\bf Novikov Betti numbers and $L_2-$Betti numbers} 
To a pair $(X, \xi\in H^1(X;\mathbb Z)),$ $X$ a compact ANR,  and $\kappa$ a field,  in addition to the familiar Betti numbers $\beta_r(X;\kappa),$ one can associate the  Novikov--Betti numbers  $\beta^N_r(X,\xi;\kappa).$  They are interesting numerical invariants of  geometric and topological relevance
\footnote {for example, 
 for  $X$ is a closed Riemannian manifold and $f:X\to \mathbb S^1$  a Morse angle valued function (i.e. all critical points non degenerated) one has   the same relation  between the numbers of critical points and the Novikov Betti numbers (Novikov inequalities) as the familiar relations between the critical points of a real valued Morse map and the standard Betti numbers (Morse inequalities)}.  
For $\kappa= \mathbb C$ the $L_2-$ Betti numbers $\beta^{L_2}_r(\tilde X),$ $\tilde X\to X$  the infinite cover determined by $\xi,$ are the same as $\beta^N_r(X,\xi;\kappa).$ 
 
Recall from \cite {N} or \cite{Pa}  that The Novikov--Betti numbers $\beta^N_r(X,\xi;\kappa)$ are defined using the infinite  cyclic cover $\tilde X\to X$ associated to $\xi$ by the equality
$$\beta^N_r(X, \xi;\kappa)= \dim _{\kappa[t^{-1},t]] }H_r(\tilde X)\otimes_{\kappa[t^{-1}, t]} \kappa[t^{-1}, t]].$$ 
Here $\kappa[t^{-1}, t]$ denotes the ring of Laurent polynomials in $t$ with coefficients in $\kappa,$   $\kappa[t^{-1}, t]]$ denotes  the field of Laurent  power series in $t$ with coefficients in $\kappa$ and $H_r(\tilde X;\kappa)$ is viewed as a $\kappa[t^{-1},t]-$module  whose $\kappa[t^{-1}, t]-$module structure determined by the action of $\mathbb Z$ on $\tilde X$ as the group of deck transformations. 
\vskip .1in  
 
Note that when $X$ is a finite simplicial complex the standard Betti numbers are {\it computer friendly} which means  computable by effective algorithms ( for instance the persistence algorithm \cite{ELZ}, \cite {ZC}, \cite {CEM}) and so are the Jordan cells in view of the algorithm presented in section \ref{S4} or of the algorithm described in \cite{BD11}. In view of the definition above, even when $X$ is a finite simplicial complex, $\tilde X$ is an infinite simplicial complex, hence the calculation of the Novikov Betti numbers via their definition is not computer friendly.

 The formula (\ref{E18})
established in \cite {BH} Theorem 7.2., 
permits however to express the Novikov--Betti numbers in terms of  {\it computer friendly} quantities and provides an alternative definition of them not based on infinite coverings.
 
 \begin {equation}\label {E18} 
 \beta^N_r(\xi;\kappa)= \beta_r(X;\kappa) - \sharp \mathcal J_r(\xi)(1) -\sharp \mathcal J_{r-1}(\xi)(1)
 \end{equation}
 In this formula  $\mathcal J_r(\xi)(1)= \{J= (\lambda, n)\in \mathcal J_r(\xi) \mid  \lambda=1\}$  and $\sharp$ denotes cardinality.

\subsection {\bf Other type of Betti numbers}

The calculation of the sets $\mathcal J_r(X;\xi)$ lead to  the calculation of the dimension of the homology with coefficients in the  local system (of one dimensional $\kappa-$vector spaces)  $\hat u \cdot \xi,\ $  $u \in\kappa\setminus 0,$  defined by the composition  

$\xymatrix { H_1(X;\mathbb Z) \ar[r]^\xi &\mathbb Z \ar[r]^{\hat u}& \kappa\setminus 0}
$ 
with $\hat u: \mathbb Z\to \kappa\setminus 0$ given by $\hat u(n)= u^n.$ 

\noindent Precisely one has the following formula  established in \cite {BH1}, Theorem 7.1,
\begin {equation} 
\dim H_r(X; \hat u \xi)= \beta^N_r(X,\xi ;\kappa) 
 + \sharp \mathcal J_r(1/u) + \sharp \mathcal J_{r-1}(u).
 \end{equation}
For $\omega\in \overline \kappa \setminus 0,$  $\mathcal J_r(\omega) $ denotes the set 
$\{J= (\lambda, n)\in \mathcal J_r(\xi) \mid  \lambda=\omega  \}.$

\subsection {\bf The Betti number of the homotopy theoretic fiber of $\xi.$} 

For $\beta_r^N(X,\xi; \kappa)=0,$  in view of the \cite {BH}, the $\kappa[t^{-1},t]-$module $H_r(\tilde X)$ is  torsion module, hence  equal to $V_r(X;\xi),$ hence  
\begin{equation}
\dim H_r(\tilde X)= \sum _{J= (\lambda_J, n_J) \mid J\in \mathcal J_r(\xi)} n_J .
\end{equation} 
If $\xi$ is represented by a fibration $f:X\to \mathbb S^1$ with compact fiber $f^{-1}(\theta)$ then $f^{-1}(\theta)$ and $\tilde X$ have the same homology, hence $\beta_r(f^{-1}(\theta);\kappa)= \dim H_r(\tilde X).$ 
By Theorem \ref {T36} 
 $\mathcal J_r(\xi)= \mathcal J_r(g)$  for $g:X\to \mathbb S^1$ any map, not necessary a fibration,  in the homotopy class  representing by $\xi.$ 
These sets can  be computed using the algorithm described in section \ref{S4}
\footnote {This formula was used by the author to calculate the Betti numbers of the Milnor fiber of some isolated singularitiesfiber for some isolated singularities and will be described in subsequent work.}.

\subsection {\bf Alexander polynomial of a knot and generalizations} 
A knot $K\subset S^3,$ $K$ a simple closed curve (i.e. homeomorphic to the oriented circle $\mathbb S^1),$ is a  locally flat embedding in the three dimensional sphere $S^3.$  Consider $X= S^3\setminus N$ where $N$ is a open tubular neighborhood of $K$ in $S^3.$ Clearly $X$ is homotopy equivalent to $S^3\setminus K.$ The Alexander dual of the generator $u\in H_1(K; \mathbb Z)$ is an integral cohomology class $\xi\in H^1(X;\mathbb Z).$
The Alexander polynomial of the knot, a fundamental invariant of the knot, is a polynomial with integral coefficients $$ a_r z^r + \cdots a_1 z+ a_0$$ with $a_0\ne 0$ and $a_n\geq 0$  is defined as the only generator of the principal ideal $(P(t))$defined by  the isomorphism  $H_r(\tilde X; \mathbb Z)\equiv \mathbb Z[t^{-1},t]/(P(t))$ where $\mathbb Z[t^{-1},t]$ denotes the ring of Laurent polynomial with coefficients in $\mathbb Z$ 

For detailed definitions and examples one recommends ( \cite{R})  \footnote { For example for the familiar figure eight knot $P(t)= t^2-3t+1,$ cf \cite {R} page 166 and the  torus knot $(4,7)$ $P(t)= t^{18} -t^{17} + t^{14}- t^{13} + t^{11} -t^9 + t^7 -t^5 +t^4 -t +1$ cf \cite{R} page178}.
 As established first by Milnor, cf \cite {Mi1}, the monic polynomial $1/a_n \cdot P(t)$ can be calculated as the characteristic polynomial of the $1-$ monodromy of $(X,\xi),$ and is exactly

\begin{equation}
\prod_{ J= (\lambda_J, n_J)\mid J\in \mathcal J_1(\xi)}  (z-\lambda_J)^{n_J}.\end{equation}
The algorithm described in Section \ref{S4} provides a new algorithm to calculate the monic Alexander polynomial \footnote{ a more in detailed discussion for the calculation of Alexanders polynomials of knots and links  including knots in higher dimension  and the role of the algorithm provided in section \ref{S4} is in preparation. 
In a similar vein important cases of  the Milnor--Turaev torsion of $(M,\xi)$, a rational function defined on $GL(n,\mathbb C),$ when regarded as the variety of rank $n$ complex representations of $Z$  can be calculated by implementable algorithms containing as part the calculation of the Jordan cell for monodromies.
 }.
 
\section{Appendices}

\subsection{\bf Appendix 1}

For the proof of Propositions \ref{P23} and \ref{P24} one needs 
the following observation.

\begin{obs} \label {O2}\ 

i).  $x\in D$   iff there exists $x_i\in V, \ i\in \mathbb Z$ with
$x_i \ R \ x_{i+1}, \ x_0=x.$ 
\vskip .1in

ii).  $y \in K_+ + K_-$  iff  there exists 
 a nonnegative integer 
$k,$  
the sequences $\ x^+_1, x^+_{2}, \cdots x^+_{k }$  all elements in V,
 and the sequence 
 $\ x^-_1, x^-_0, x^-_{-1}, \cdots x^-_{-k}$  all elements in V,  
such that:  
\begin{enumerate}

\item$y= x_1^+ + x_1^-,$ 
\item $x^+_1\ R\ \ x^+_2\ R\cdots x^+_k\ R\ \ 0, $
\item $0\ R\ x^-_{-k}\ R\  x^-_{-(k-1)}\ R\ \cdots x^-_{0} \ R\  \ x^-_1.$
\end{enumerate}
\end{obs}

{\bf Proof of Proposition \ref{P23}} (cf \cite {BH})

To establish item 1.  one uses Lemma  \ref{L22} (2) applied to  the relation $R_{reg}.$
Clearly in view of  the surjectivity of $\pi :D\to  V_{reg}$ and Observation \ref{O2} (i) one has
$\dom R_{reg} =V_\reg,$ so it remains to check  that $\ker (R_\reg)=0.$ 

To verify this  we start with 
$x\in D$ s.t  $x R   y, \ \ y\in (K_+ +K_-)$ and want to check that $x\in D\cap (K_+   + K_-).$

One produces the elements  $x_i \in V,  i\in \mathbb Z,$  $\ x^+_1, x^+_{2}, \cdots x^+_{k }\in V$  and the elements $\ x^-_1, x^-_0, x^-_{-1}, \cdots x^-_{-k}\in V$  as stated in Observation \ref{O2} (ii)  and with the properties  stated.  One  observes that:

\begin{enumerate}
\item $x^-_0\in K_-,$
\item 
$(x-x^-_0)\in D \cap K_+$  since 
$$ \cdots R x_{(-k-1)} R (x_{-k}- x^-_{-k}) R \cdots R(x_0- x^-_0) R\ ((y-x_1^-)= x^+_1)\ R x^+_2 \cdots R x^+_k  R \  0$$
and  therefore 
$x^-_0= -x+ x- x_0^- \in D,$ hence 
\item $x^-_0 \in D\cap K_-,$  
\end{enumerate}
 
 Combining (2.) and (3.) above  one obtains  
$x= x-x_0 +x_0\in (D\cap K_+)+ (D\cap K_-) \subseteq D\cap (K_+ + K_-).$ 
\vskip .1in 
Items  2. 3. and  4. (in Proposition \ref{P23}) are straightforward. 

To verify item 5. it suffices to check the equality for $k=2.$  which 
can be concluded  in view of Observation \ref{O2} (i). 
\vskip .1in 

{\bf Proof of Proposition \ref{P24}} 

Item 1. follows by observing that $D$ and $D\cap(K^++K^-)$ for both $R(\alpha, \beta)$ and $R(\alpha', \beta')$ are actually the same.

To check this 
consider the sequences 

\scriptsize{
\begin{equation*}
\xymatrix {
&...v_{-1}\ar[r]^\alpha&w_{0}&v_0=(v^-_0 +v_0^+)
\ar[l]_\beta\ar[r] &w_1&v_1\ar[l]\ar[r]&w_2\cdots 
\\
&v_0^+\ar[r]&w^+_1&v^+_1\ar[l]\ar[r]&w^+_2&\cdots\ar[l] 
w^+_{k+1}&0\ar[l]
\\
&0\ar[r]&w^-_{-(k)}&v^-_{-k}\ar[r]\ar[l]&w^-_{-(k-1)}&\cdots
v_{-1}\ar[l]\ar[r]&w^-_{0}
}
 \end{equation*}}
 \normalsize

Indeed, by Observation \ref {O2} (i),  $v_0\in D$ implies  the existence of the first sequence, which implies  
 that $v_i\in V'$ and $w_i\in W',$  which guarantee that $D=D'.$ 

If $v_0\in D\cap (K_+ + K_-)$ all  three  sequences above exist, which imply that 
 that $v_0- v^-_0 = v_o^+\in D\cap K'_+ \subseteq D'\cap (K'_+ K'_+).$ Similarly  $v_0- v^+_0= v_0^-\in D'\cap K'_-\subseteq D'\cap (K'_+ + K'_-),$ and therefore $v_0= v- v^-_0 + v - v^+_0= v_0\in D'\cap ((K'_+ + K'_-).$

\vskip .1in
To check item  2. observe that
the diagram (2) (in Section \ref{S1}) induces the linear map $\pi : D/ D\cap (K_+ + K_-)\to D'/ D'\cap (K'_ - + K'_+).$ This map is obviously surjective since both pairs $\alpha,  \beta$ and $\alpha' , \beta'$ being surjective make $V=D$ and $V'=D'.$ 
To check that is injective we will verify that   ${p'}^{-1}(K'_\pm)\subset K_\pm.$ 

For this purpose consider  diagram (2)   with $\alpha' $ and $\beta'$ as specified by hypotheses.

\begin{lemma} \label {L52} 
If $w\in W, w'\in W', v'\in V'$ such that $p(w)= w'$ and $\beta'(v')= w'$ then there exists $v\in V$ such that $\beta(v)=w$ and $p'(v)=v'.$ 
\end{lemma}

\proof We first choose $\underline v$ with the property $p'(\underline v)=v',$ observe that $p(w-\beta(\underline v))=0,$ hence 
in view of the definition of the diagram (2) $w-\beta(\underline v)= \beta(u), u\in \ker \alpha$ and correct $\overline v$ to $v$  by taking $v= \underline v -u.$ 

q.e.d
\vskip .1in 

With Lemma \ref{L52} established observe that given a sequence 
 $v'_0, v'_1, \cdots v'_k \in V'$  and $v_0 \in V$ with the property that 
\begin{equation}\label {E11}
\begin{aligned}
\alpha' (v'_{i-1})=& \beta'(v'_{i}),\  1 \leq i \leq k  \\
p(v_0)= &v'_0
\end{aligned}
\end{equation}

one can produce 
$v_1, v_2, \cdots v_k \in V$
such that 

\begin{equation}\label {E12}
\begin{aligned}
\alpha(v_{i-1})= &\beta(v_{i}) \\
p(v_i)= &v'_i .
\end{aligned}
\end{equation}

Indeed suppose inductively that $v_1, v_2,\cdots  v_i, \ i\leq r$ satisfying  properties (\ref{E12})  are produced.  Apply Lemma \ref{L52} to $w=\alpha(v_i),
w'= \alpha'(v'_i)$ and $v'= v'_{r+1}$ and obtain  $v_{r+1}.$
\vskip .1in
To conclude ${p'}^{-1}(K'_+) \subset K_+$  one chooses the sequence $\{v'_i\}$ to have (for some  $k$) $\alpha(v'_k)=0,$  which means that $v'_0 \in K'_+.$ Then $v_k$ constructed as above is in $ \ker \alpha$ which means that $v_0\in K_+.$ 

To conclude ${p'}^{-1}(K'_-) \subset K_-$ one chooses a sequence  $\{v'_i\}$  to have (for some $k$)  $v'= v'_k \in K'_-$  and $v'_0=0$  Then for for $v_0=0$ one construct the sequence $v_1, v_2, \cdots v_k \in V.$ Then
$v_k\in K_-,$  hence ${p'}^{-1}(K'_-)\subset  K_- .$    

q.e.d.

\subsection{\bf Appendix 2.}

Recall that: 

 The Hilbert cube $Q$ is the infinite product $Q= \prod_{i\in \mathbb Z_{\geq 0}} I_i= I^{\infty}$  with $I_i=I=  [0,1].$ 
The  topology of $Q$ is  given by the metric $d( u ,v)= \sum _i |u_i- v_i|/ 2^i$ with $ u= \{u_i\in I, i\in \mathbb Z_{\geq 0}\}$  and  $ v= \{v_i\in I, i\in \mathbb Z_{\geq 0}\}.$

A compact Hilbert cube manifold is a compact Hausdorff space locally  homeomorphic to the Hilbert cube.

The space $Q$ is a compact ANR and so is $X\times Q$ when $X$  is a compact ANR.
\vskip .1in 

\noindent For any $n,$ positive integer, write $Q= I^n\times Q'_n$ and denote by
$\pi_n:Q\to  I^n$ the first factor projection and by
 $\pi_n^X: X\times Q\to X\times I^n$ the product $\pi_n^X= id_X\times \pi_n.$

\noindent  For $F:X\times Q\to \mathbb R$ let $F_n$ be the restriction of $F$ to $X\times I^n$ and  $\overline F_n$
the composition $\overline F_n: =F_n\cdot \pi^X_n.$ 

\noindent For $f:X\to \mathbb R$ denote by $\overline f := f\cdot \pi_X$ where $\pi_X: X\times Q\to X$  is the canonical projection on $X.$ 
Note that  

\begin{obs}\label {OA}\
\begin {enumerate}
\item If $f:X\to \mathbb R$ is a tame map so is $\overline f.$  
\item  The sequence of maps $\overline F_n$ is  uniformly convergent to  the map $F.$ 
\end{enumerate}
\end{obs}

The following are two  results about Hilbert cube manifolds whose proof can be found in \cite {CH}.

\begin{theorem} \label {T52}\
\begin{enumerate} 
\item (R Edwards) If $X$ is a compact ANR then $X\times Q$ is a Hilbert cube manifold.
\item (T Chapman) If $\omega:X\to Y$ is a homotopy of equivalence  between two finite simplicial complexes with Whitehead torsion $\tau(\omega)=0$ then   there exists a homeomorphism $\omega': X\times Q\to Y\times Q$ such that 
$\omega'$ and $\omega\times id_Q$ are homotopic.   
\end{enumerate}
\end{theorem}

Recall  (for the non expert reader)  that for a homotopy equivalence $f:X\to Y$ between two finite simplicial complexes  one can associate an element  $\tau (f)\in Wh(\pi_1(X,x)) $  which measures the obstruction to $f$ to be a "simple homotopy equivalence" in the sense of J.H. Whitehead. Here  $Wh(\Gamma)$ denotes the Whitehead group of $\Gamma,$  which is an abelian group  associated with a discrete group $\Gamma, $ cf \cite{Mi}. It is also known \cite {Mi} that if $K$ is a finite cell complex (actually a compact ANR)  with $\chi(K)=0$ then $\tau  (f\times Id_K)=0$ in particular  $\tau  (f\times Id_{\mathbb S^1})=0.$ In view of theorem above Chapman has extended $\tau(f)$ to a homotopy equivalence of compact Hilbert manifolds. 
\vskip .1in

\noindent { \bf Proof of Stabilization theorem:}

Items 1. and 2. in Stabilization theorem follow from item 1. respectively item 2. combined with item 3. in Theorem \ref{T52}.   
\vskip .1in

One has also the following result  whose proof was provided by S. Ferry:
\begin{proposition} A compact Hilbert cube manifold is a '' good  ANR ''.
\end{proposition}

\proof  \ Let $M$ be a Hilbert cube manifold and $F:M\to \mathbb R$ a continuous map. We want to show that for $\epsilon >0$ one can produce a tame map $P : M\to \mathbb R$ such that  $|F( u)- P( u)|<\epsilon$ for  any   $ u\in M.$
For this purpose write $M= K\times Q,$ $K$ a finite simplicial complex, cf \cite{CH} section 11. 

It suffices to produce an $n$ and a simplicial map $p: K\times I^n\to \mathbb R$ such that $|F- p\cdot \pi^K_n|<\epsilon.$ 

The continuity of $F$ and the compacity of $M$ insure the existence of $\delta>0$ such that $| u- v|<\delta$ implies $|F( u)- F( v)|<\epsilon/2,\ \ $  $u,  v \in M= K\times Q.$

Choose $n$ such that  $| u-(\pi^K_n( u),0)|<\delta, \   u\in K\times Q $ (here $(\pi^K_n( u), 0)\in (K\times I^n)\times Q'_n= K\times Q$).

Choose $p:K\times I^n\to \mathbb R$ a simplicial map (with respect to a convenient subdivision) with $|p(x)- F_n(x)|<\epsilon/2, \ x\in K\times I^n,$ and take $P= p\cdot \pi^K_n.$ Since $p$ is tame so is $P.$ 

Then $|F( u)- p\cdot \pi^X_n( u)|\leq |F( u)- F_n\cdot \pi^K_n( u)|+ |F_n\cdot \pi^X( u)- p\cdot \pi^K_n( u)| < \epsilon.$
q.e.d.

\subsection {\bf Appendix 3}

Recall from \cite {BD11} sections 4 and 5 
the following notations: 

-- The oriented graph $G_{2m}$ has vertices $x_1, x_2, \cdots x_{2m}$ and the oriented edges 
$a_i: x_{2i-1}\to x_{2i}, \  b_i:x_{2i+1}\to x_{2i} $ with $x_{2m+1}= x_1, i=1,\cdots m.$  

-- A $G_{2m}-$representation  $\rho$ is given by a collection of 
linear maps $\alpha_i: V_{2i-1} \to V_{2i}, \beta_i: V_{2i+1}\to V_{2i}$ with $V_i$ vector space corresponding to the vertex $x_i,$  and the linear map $\alpha_i$ resp. $\beta_i$ corresponding to the arrow $a_i$ resp. $b_i.$ 

--For  $f:X\to \mathbb S^1$  a tame map in the sense of \cite{BD11} with $m$ critical angles $0 < s_1< s_2< \cdots s_m\leq 2\pi$ and $t_1,  t_2,  \cdots \ t_m $  regular values such that $0< t_1< s_1 < t_2 \cdots s_{m-1}<t_m <s_m$  one associate the $G_{2m}-$representation $\rho_r$ with  $V_{2i-1}= H_r(\tilde f^{-1}(t_i)),$   $V_{2i} = H_r(f^{-1}(s_i))$  and $\alpha^r_i, \ \beta^r_i$ the linear maps induced in homology by the continuous maps $a_i$ and $b_i,$ considered in \cite{BD11} Section 5. 
\vskip .1in

Let   $\tilde f:\tilde X\to \mathbb R$ be the canonical  infinite cyclic cover of the tame map $f:X\to \mathbb S^1.$ 
Put $t_{m+1}= t_1 +2\pi$ and observe that 
$V_{2i}=  
H_r({\tilde f}^{-1}(t_t, t_{i+1}))$ and the relation $R^{\tilde f,\tilde f}_{t_i, t_{i+1}}(r)$ is exactly 
$$R^{\tilde f,\tilde f}_{t_i, t_{i+1}}(r)= R(\alpha^r_i, \beta^r_i)= R(\beta^r_i)^\dag \cdot R(\alpha^r_i).$$
Clearly  the composition $$R^{\tilde f, \tilde f}_{t_m, t_{m+1}}(r)\cdot R^{\tilde f, \tilde f}_{t_{m-1}, t_m}(r)\cdots R^{\tilde f, \tilde f}_{t_2, t_3}(r)\cdot R^{\tilde f, \tilde f}_{t_1, t_2}(r)$$  identifies to $R^f_{t_1}(r).$

\vskip .1in 

To a $G_{2m}-$ 
representation $\rho$ one associates the linear relation $R(\rho)\colon V_1\rightsquigarrow V_1=  R^\dag(\beta_m)\cdot R(\alpha_m)\cdots R^\dag(\beta_1)\cdot R(\alpha_1)$ and one denotes by $\mathbb J  (\rho):= \mathcal J([R(\rho)_\reg]).$
Clearly  one has:
\begin{obs} \label {OA3}\
\begin{enumerate}
\item $R(\rho\oplus \rho')= R(\rho)\oplus R(\rho')$ and therefore $\mathbb J(\rho)\sqcup \mathbb J(\rho'),$
\item $\mathbb J(\rho^I)= \emptyset,$
\item $\mathbb J(\rho^I(\lambda, k))= \{ (\lambda,k) \}.$
\end{enumerate}
\end{obs}

The Jordan cells $(R^f_{t_1})_\reg (r)$ are the Jordan cells of the $r-$monodromy  $T^{(X, \xi_f)}(r)$ and then, by Observation \ref{OA3}, they are  are the Jordan cells of the representation $\rho_r$ defined in \cite{BD11}.

\end{document}